\documentclass[10pt]{amsart}
\usepackage{epsfig}
\usepackage{amsmath}
\usepackage{amssymb}
\usepackage{xr}
\usepackage{amscd}

\newcommand{\ord}{\mathcal{O}}
\newcommand{\re}{\text{Re}}

\newtheorem{lem}{LEMMA}[section]
\newtheorem{theo}[lem]{THEOREM}

\newtheorem{definition}[lem]{DEFINITION}
\newtheorem{remark}[lem]{Remark}

\begin{document}
\author{Daniel Jupiter, Krastio Lilov}
\title[Invariant Fatou Components]{Invariant nonrecurrent Fatou
components of automorphisms of $\mathbb{C}^2$}{

\begin{abstract}
Let $\Omega$ be an invariant nonrecurrent Fatou component
associated with the automorphism
$F:\mathbb{C}^2\rightarrow\mathbb{C}^2$. Assume that all of the
limit maps of $\{F^n|_{\Omega}\}$ are constant. We prove the
following theorem. If there is more than one such limit map then
there are uncountably many. The images of these limit maps form a
closed set in the boundary of $\Omega$ containing no
isolated points. Additionally there cannot be more than one limit
map unless the derivative of $F$ along a specific subset of the
curve of fixed points of $F$ has eigenvalues $1$ and
$e^{i2\pi\theta}$, with $\theta$ non-Diophantine.

We also examine the case where the limit maps are not all constant.
The image of a nonconstant limit map is an immersed
variety in the boundary of $\Omega$. We show that any two such
immersed varieties intersect either trivially or in a set that is
open in their intrinsic topologies.

We present some examples of maps with invariant nonrecurrent Fatou
components.
\end{abstract}

\maketitle

\section{Introduction}
The Fatou components for rational self maps of $\overline{\mathbb{C}}$ are
entirely classified: see for example
\cite{carleson-gamelin:dynamics}. Forn{\ae}ss and Sibony
\cite{fornaess-sibony:fatou2} have examined recurrent Fatou
components of holomorphic self maps of $\mathbb{P}^2$ of degree at
least $2$. Ueda \cite{ueda:fatou} has also made contributions in this
direction. Forn{\ae}ss and Sibony have
also \cite{fornaess-sibony:fatou} studied recurrent Fatou
components for generic maximal rank $k$ holomorphic self maps of
$\mathbb{C}^k$. Bedford and Smillie \cite{bedford-smillie:1} and
\cite{bedford-smillie:2} have investigated Fatou components of
H{\'e}non maps. In this paper we consider invariant nonrecurrent
Fatou components, $\Omega$, for an automorphism, $F$, of
$\mathbb{C}^2$.

We consider first the case where all the limit maps are rank $0$.
Here we show that generically there cannot be more than one such
limit map. Next we examine the case where there are rank $1$ limit
maps. We show that the images of such maps generically do not
intersect. We show next that a large class of polynomial
automorphisms of $\mathbb{C}^2$ do not have invariant nonrecurrent
Fatou components. Finally we provide some examples of maps with
invariant nonrecurrent Fatou components.

In Section \ref{rank 0} we examine the case where all limit maps
are rank $0$. We say that such a Fatou component satisfies
Property $0$. We let $J$ be the set of fixed points of $F$,
$\Sigma$ be the images of the limit maps of $\{F^n\}_{n=0}^{\infty}$. We
construct an $F$ invariant curve which lies entirely in $\Omega$.
If certain eigenvalues of $F'$ along $\Sigma$ satisfy the
Diophantine condition (\ref{Siegel condition}), defined on page
\pageref{Siegel condition}, we construct continuously varying
families of invariant manifolds transversal to $J$. Using the
invariant curve and the invariant manifolds which we have
constructed, we prove the following theorem.
\begin{theo}\label{main theorem}
Let $\Omega$ be an invariant, nonrecurrent Fatou component
satisfying Property $0$. If $\Sigma$ contains more than one point,
then the following hold:
\begin{enumerate}
\item $\Sigma$ is uncountable and has no isolated points. \item
There is a one dimensional subvariety
$V\subset\{(z,\,w)\in\mathbb{C}^2\,;\,F(z,\,w)=(z,\,w)\}$ such
that $\Sigma\subset V\cap\partial\Omega$. \item There exists an
$\alpha\in\mathbb{C}$, with $\alpha=e^{i2\pi\theta}$, with
$\theta$ non Diophantine, such that the eigenvalues of $F'(z,\,w)$
are $\{1,\,\alpha\}$ for all $(z,\,w)\in\Sigma$.
\end{enumerate}
\end{theo}

In Section \ref{rank 1} we consider the case where the limit maps
may be either rank $0$ or rank $1$. In Section
\ref{section:composing limit maps} we find a natural way to extend
a rank $1$ limit map, $h$, to its image, allowing us to examine
the action of iterates of $F$ on this image. We conclude that the
family $\{F^n\}$ is normal on $h(\Omega)$, if we consider $F^n$ as
a map from $h(\Omega)$ to $\mathbb{C}^2$. Using the fact that the
images of $\Omega$ under rank $1$ limit maps are immersed
varieties, in Section \ref{section:limit varieties} we prove the
following theorem.
\begin{theo}\label{varieties}
Let $h(\Omega)$ and $g(\Omega)$ be two distinct rank $1$ limit
varieties. Then $h(\Omega)\cap g(\Omega)$ is either empty or an
open set, when considered as a subset of $h(\Omega)$ or
$g(\Omega)$.
\end{theo}

In Section \ref{polynomial} we examine the question of whether
polynomial automorphisms can have periodic nonrecurrent Fatou
components, and prove the following theorem.
\begin{theo}\label{theorem:polynomial}
If $F$ is a polynomial automorphism of $\mathbb{C}^2$ then it
cannot have an invariant nonrecurrent Fatou component on which it
has more than one rank $0$ limit map.

If $F$ is a polynomial automorphism of $\mathbb{C}^2$ with an
invariant nonrecurrent Fatou component on which it has one rank
$0$ limit, then it is a H{\'e}non map.

If $F$ is a polynomial automorphism of $\mathbb{C}^2$ with an
invariant nonrecurrent Fatou component on which it has a rank $1$
limit, then it is a nonhyperbolic H{\'e}non map.
\end{theo}

In Section \ref{examples} we present some examples of such maps.
We present an example of an automorphism which has an invariant
nonrecurrent Fatou component with exactly one rank $0$ limit map.
We give two examples of automorphisms with rank $1$ limits. The
first has precisely one rank $1$ limit, with image the $w$ axis.
The second rank $1$ example has multiple rank $1$ limit maps, all
with image the $w$ axis.

 The authors would like to thank
Professors Eric Bedford, John Erik Forn{\ae}ss and Berit
Stens{\o}nes for helpful advice and comments in writing this
paper.

\section{The Rank $0$ Case} \label{rank 0}

Let $F:\mathbb{C}^2\rightarrow\mathbb{C}^2$ be an automorphism.
Let the Fatou set, $\mathcal{F}$, denote the set of points where
the forward iterates $\{F^n\,;\;n\geq 0\}$ are locally normal:
every subsequence $\{F^{n_j}\}$ has a further subsequence that
converges uniformly on compact subsets of $\mathcal{F}$. Let
$\Omega$ denote a connected component of $\mathcal{F}$, i.e. a
Fatou component. A Fatou component $\Omega$ is said to be
invariant if $F(\Omega)=\Omega$. We say that $\Omega$ is recurrent
if there is a point, $z_0\in\Omega$, and a sequence of integers
$\{n_j\}$, such that $F^{n_j}(z_0)\rightarrow z\in\Omega$.

If $\Omega$ is a Fatou component, let $\Sigma$ denote the set of
maps, $h:\Omega\rightarrow\mathbb{C}^2$, which are obtained as
normal limits of subsequences of $\{F^n|_{\Omega}\}$. If $\Omega$
is invariant but not recurrent, then every $h\in\Omega$ maps
$\Omega$ to $\partial\Omega$. In this section we consider the case
where each $h\in\Sigma$ has rank zero, and is thus constant. We
say in this case that $\Omega$ satisfies Property $0$. We identify
$h$ with the point $h(\Omega)\in\partial\Omega$.

The remainder of Section \ref{rank 0} consists of a proof of
Theorem \ref{main theorem}.

\begin{lem}\label{rank 0 theorem}
Let $F$, $\Omega$  and $\Sigma$ be as in the statement of Theorem
\ref{main theorem}. Then we have the following.
\begin{enumerate}
\item $\Sigma$ contains uncountably many elements. \item $\Sigma$
is a closed set containing no isolated points. \item Every point
in $\Sigma$ is a neutral fixed point.
\end{enumerate}
\end{lem}

\begin{proof}
The Lemma is proved by constructing an $F$ invariant curve lying
in $\Omega$, which connects the points $\{F^n(p)\}$, for some
$p\in\Omega$. Given the images of any two limit maps, the curve
travels back and forth between these images, and clusters at the
boundary of $\Omega$. These cluster points are the uncountably
many limit points mentioned in the Lemma. Using the local normal
forms described by Ueda \cite{ueda:local} we show that each of
these points is a neutral fixed point.
\subsection{An Invariant Curve}\label{invariant curve}
We construct an $F$ invariant curve $\gamma$ which we use to show
that if we have two limit maps in $\Sigma$ then in fact we have
uncountably many. Let $p\in\Omega$. Let
$\gamma_0:[0,1]\rightarrow\mathbb{C}^2$ be any smooth curve from
$p$ to $F(p)$ which is entirely contained in $\Omega$. Let
$\gamma_n:[n,n+1]\rightarrow\mathbb{C}^2$ be $F^n(\gamma_0(t-n))$.
Finally, let $\gamma:[0,\infty)\rightarrow\mathbb{C}^2$ be given
by
\[\gamma(t)=\gamma_n(t) \text{ when } t\in[n,n+1].\]
Clearly $\gamma$ is $F$ invariant.

Let $q=\lim_{j\rightarrow\infty}\left(F|_{\Omega}\right)^{m_j}$
and $q'=\lim_{k\rightarrow\infty}F^{n_k}$ be two distinct points
in $\Sigma$, and let $B^2(q,\epsilon)$ and $B^2(q',\epsilon)$ be
two nonintersecting balls around $q$ and $q'$.\\ Notice that
$\gamma$ is connected. Also note that for large enough $j$ we have
$F^{m_j}(\gamma_0)\subset B^2(q,\epsilon)$ and
$F^{n_k}(\gamma_0)\subset B^2(q',\epsilon)$. By the connectivity
of $\gamma$ there must be points of $\gamma$ in $\partial
B^2(q,\epsilon)$. We make this more precise.\\ Given any $t_0 \in
[0,1]$ there are an $m'_0\in\left\{m_j\right\}$ and an
$n'_0\in\left\{n_k\right\}$ such that
\begin{enumerate}
\item $m'_0<n'_0$, \item
$F^{m'_0}(\gamma_0(t_0))=\gamma_{m'_0}(t_0+m'_0)\in
B^2(q,\epsilon)$, and \item
$F^{n'_0}(\gamma_0(t_0))=\gamma_{n'_0}(t_0+n'_0)\in
B^2(q',\epsilon)$.
\end{enumerate}
Then there is a $t_1$ such that
\begin{enumerate}
\item $m'_0<t_1<n'_0$, and
\item $\gamma(t_1)\in \partial B^2(q,\epsilon)$.
\end{enumerate}
Repeating the above procedure, but choosing $m'_1$, $n'_1$ both
bigger than $n'_0$, we produce $m'_1<t_2<n'_1$ such that
$\gamma(t_2)\in \partial B^2(q,\epsilon)$ We repeat again,
producing the sequences $\left\{t_i\right\}_{i=0}^{\infty}$,
$\left\{m'_i\right\}_{i=0}^{\infty}$,
$\left\{n'_i\right\}_{i=0}^{\infty}$, and
$\left\{\gamma(t_i)\right\}_{i=0}^{\infty}$. We note that for each
$t_i$ we have
\[\gamma(t_i)=\gamma_{l_i}(\zeta_i+l_i)=F^{l_i}(\gamma_0(\zeta_i)),\]
for some integer $l_i$ and some $\zeta_i\in [0,1]$. Passing to
subsequences we can assume that
\begin{enumerate}
\item $\zeta_i\rightarrow\zeta\in [0,1]$,
\item $\gamma(t_i)\rightarrow\eta_{\epsilon}\in \partial B^2(q,\epsilon)$, and
\item $\left(F|_{\Omega}\right)^{l_i}$ converges.
\end{enumerate}
Then
\begin{align*}
\lim_{i\rightarrow\infty} F^{l_i}(\gamma_0(\zeta))
    &= \lim_{i\rightarrow\infty} F^{l_i}(\gamma_0(\zeta_i))\\
    &= \lim_{i\rightarrow\infty}\gamma(t_i)\\ &= \eta_{\epsilon}.
\end{align*}
We know that $\eta_{\epsilon}\in \partial\Omega$: if it were not,
$\Omega$ would be recurrent.

Given two limit maps $q$ and $q'$ we have produced a third,
$\eta_{\epsilon}$. We notice, however, that for each $\epsilon$
suitably small we find a different $\eta_{\epsilon}$. We have thus
shown that given two limit maps there are in fact uncountably
many.

\subsection{The Structure of $F$ on $\Sigma$}\label{structure of
sigma} By the construction of the previous section, given a limit
map $q$, there are other limit maps arbitrarily close to $q$. In
other words, $\Sigma$ contains no isolated points.

We claim that $\Sigma$ is a closed set. To prove this assume there
is $q\in\overline{\Sigma}\backslash\Sigma$. Then, given
$\epsilon>0$, there is a point $q_{\epsilon}\in\Sigma$ such that
\[|q-q_{\epsilon}|<\epsilon/2.\]
Since $q_{\epsilon}\in\Sigma$, we can find $n_{\epsilon}$ so that
for a fixed $p\in\Omega$ we have
\[|F^{n_{\epsilon}}(p)-q_{\epsilon}|<\epsilon/2.\]
But then
\[|F^{n_{\epsilon}}(p)-q|<\epsilon.\]
Passing to a convergent subsequence of
$\left\{\left(F|_{\Omega}\right)^{n_{\epsilon}}\right\}$, say
$\left\{\left(F|_{\Omega}\right)^{n_j}\right\}$, we see that
$\left(F|_{\Omega}\right)^{n_j}\rightarrow q$.

Any point
$q=\lim_{j\rightarrow\infty}\left(F|_{\Omega}\right)^{n_j}\in\Sigma$
is fixed:
\begin{align*}
F(q) &= F\lim_{j\rightarrow\infty}F^{n_j}(p) \\
     &= \lim_{j\rightarrow\infty}F^{n_j}(F(p)) \\
     &= q.
\end{align*}

We denote the set of fixed points of $F$ by $J$. $J$ is an
analytic variety in $\mathbb{C}^2$, so it is either an open set,
or consists of a union of isolated points and one complex
dimensional curves. The $F$ we are examining cannot be the
identity, so $J$ cannot be open. We have also seen that $\Sigma$
contains no isolated points. Thus we see that any point in
$\Sigma$ lies on a one dimensional complex curve in $J$. According
to Ueda \cite{ueda:local} Section 3 we can thus find local
coordinates around any point $q\in\Sigma$ in which $F$ takes the
form
\begin{equation}\label{normal form 1}
F(x,y)=(x+g(x,y)h(x,y),y+g(x,y)k(x,y)),
\end{equation}
where $g(x,y)$ is a defining equation of $J$. Singular points of
$J$ are discrete, so we can find $q\in\Sigma$ at which $J$ is not
singular. Then, again according to Ueda \cite{ueda:local}, we may
choose local coordinates in which $J$ is the $x$ axis:
\[F(x,y)=(x+yh(x,y),y(1+k(x,y))).\]
We note that along the $x$ axis we have
\[DF(x,0)=\begin{pmatrix}1 & h(x,0) \\ 0 & 1+k(x,0)\end{pmatrix}.\]
The eigenvalues of $DF(x,0)$ are $1$ and $1+k(x,0)$. There are
three possibilities:
\begin{enumerate}
\item $q$ is semi repulsive: $|1+k(x,0)|>1$,
\item $q$ is semi attractive: $|1+k(x,0)|<1$, or
\item $q$ is neutral: $|1+k(x,0)|=1$.
\end{enumerate}
To show that the former two cases are not possible, we note that
Nishimura \cite{nishimura:automorphisms} has shown that in the semi
repulsive
(resp. semi attractive) case $F$ can be written, in suitable
coordinates, as
\begin{equation}\label{normal form 2}
F(x,y)=(x,b(x)y),
\end{equation}
with $|b(x)|>1$ (resp. $<1$).

To show that $q$ is not semi repulsive, assume that $|b(x)|>1$ in
some suitably small neighbourhood, $V$, of $q=(0,0)$. Assume as
well, by shrinking $V$ if needed, that $|b|$ is bounded above on
$V$ and that $V$ is a polydisk of polyradius $\epsilon$. Fixing a
point $(x,y)\in V\cap\Omega$, with $y\neq 0$, we know that
\[(x_{n_j},y_{n_j})=F^{n_j}(x,y)\rightarrow q,\]
for some subsequence $\left\{n_j\right\}$ of integers. So for all
$j$ suitably large we have that
\[|y_{n_j}|<\frac{\epsilon}{2\max_{(x,y)\in\overline{V}}|b(x)|^2}.\]
Carefully choosing $m_j$, since $|b(x)|>1$ uniformly on $V$ and
$y_{n_j}\neq 0$, we can
 arrange that
\[(x_{n_j+m_j},y_{n_j+m_j})=(x_{n_j},(b(x_{n_j}))^{m_j}y_{n_j})\in
V,\] and that
\[\frac{\epsilon}{|b(x_{n_j})|^2}\leq |y^{n_j+m_j}| \leq
\frac{\epsilon}{|b(x_{n_j})|}.\] This gives us a new sequence,
$\left\{n_j+m_j\right\}$, where $F^{n_j+m_j}(x,y)$ lie in
\[\left\{(x,y)\in V\mid\frac{\epsilon}{|b(x)|^2}\leq |y| \leq
\frac{\epsilon}{|b(x)|}\right\}.\] This set is a compact, so
$\left\{F^{n_j+m_j}(x,y)\right\}$ has a limit point, $q$, in
$\overline{V}$ which is not in the $x$ axis. If  $q \in \Omega$
then $\Omega$ is a recurrent Fatou component, which we have
assumed it is not. If $q$ is in the boundary of $\Omega$, then
since $\Omega$ satisfies Property $0$ it is the image of a limit
map of $\left\{\left(F|_{\Omega}\right)^n\right\}_{n=1}^{\infty}$.
By above results $q$ is a fixed point of $F$. But, if $V$ is small
enough, the only fixed points in $\overline{V}$ are on the $x$
axis. We see that our fixed points cannot be semi repulsive.

Assume by way of contradiction that $q$ is a semi attractive fixed
point. It is clear by the form of Equation (\ref{normal form 2})
that in a small neighbourhood of $q$ we have convergence of the
iterates of $F$ to a map whose image contains an open set in the
$x$ axis. But we assumed that all limits of $F$ were constant
maps, so this situation is not possible. Another way of looking at
this is that the normal form also shows that we have normality in
a neighbourhood of our fixed point. But our fixed point is
supposed to be in the boundary of the Fatou component.

Note that we have not accounted for points in $\Sigma$ which are
singular points of $J$. However, both by continuity and by
Equation (\ref{normal form 1}) we see that these points are also
neutral.

This completes the proof of the lemma.
\end{proof}

\begin{remark}
The following remark holds unless the eigenvalues of $F'$ are both
constant along the curve of fixed points.

We note that at smooth points of $J$ we have that $1+k(x,0)$ is
holomorphic: denote the eigenvalue $1+k(x,0)$ by $\lambda(z)$.
Then $\lambda(z)$ is a root of the polynomial
$P(\lambda,z)=\det(F'(z)-\lambda I)$. This polynomial has
holomorphic coefficients, and the root $1$. Thus $P$ can be
written as $(\lambda-1)(\lambda-\lambda(z))$. We thus see that
$\lambda(z)$ is the constant term of $P$, and thus holomorphic.

Now we recall that Ueda tells us that the eigenvalue in the
direction of $J$ is constant $1$, and that the eigenvalue in the
transverse direction to $J$ is $1+k(x,0)$. On $\Sigma$ we know
that $|1+k(x,0)|=1$, and thus $1+k(x,0)$ varies real analytically,
away from singularities of $J$.

We notice that the above considerations show that $\Sigma$ is
contained in a locally finite union of local real analytic curves,
with discrete singularities; precisely the curves where one
eigenvalue of $DF$ is exactly $1$ and the second eigenvalue of
$DF$ has modulus $1$.
\end{remark}

We consider one-parameter families of (local) holomorphic
diffeomorphisms of $\mathbb{C}^2$ and study the parametric
dependence of local invariant manifolds. We study a case where the
maps are not necessarily hyperbolic at the fixed point. This will
be the key tool in the proof of Theorem \ref{main theorem}.

We denote coordinates on $\mathbb{C}^2$ by $z=(z_1,\,z_2)$. We
shall need the following Diophantine condition
\begin{equation}\label{Siegel condition}
    |\lambda^k-1|>ck^{-N},\ k=1,2,\ldots,
\end{equation}\\
where $\lambda=e^{2\pi i \theta}$.

\begin{theo}\label{theorem:invariant manifolds}
Let $F=F_r:(\mathbb{C}^2,\,0)\rightarrow (\mathbb{C}^2,\,0)$ be a
family of local holomorphic diffeomorphisms with the following
properties:

\begin{enumerate}
\item $F(z)=(\lambda z_1 +O(|z|^2),\, z_2+O(|z|^2))$,

\item $F$ depends holomorphically on $\lambda$.
\end{enumerate}

We restrict our attention to $\lambda=r e^{2\pi i \theta_0}$,
$1-\delta<r<1+\delta$ for $\delta$ small to be chosen later, and
$0<\theta_0<1$ is a fixed irrational satisfying the Diophantine
condition (\ref{Siegel condition}).

Then:
\begin{enumerate}

\item For any $r$ in the above range there exists a local invariant manifold
$\psi=\psi_r: \Delta(0,\,\rho)\rightarrow \mathbb{C}^2\, \ F\circ
\psi_r(w)=\psi_r(\lambda w),\ \psi_r(0)=(0,\,0)\
\psi_r'(0)=(1,\,0)$, for a fixed $\rho$ independent of $r$.
\item There exist uniform bounds
$|\psi_r|_{L_\infty(D(0,\rho))}<C$ independent of $r$.
\item The family $\psi_r$ is $C^1$-smooth in $r$ in the compact-open topology of
$\mathcal{O}(\Delta,\,\mathbb{C}^2)$.
\end{enumerate}
\end{theo}

\begin{proof}
The proof of the Theorem is a parameterized version of Siegel's
linearization Theorem, and uses ideas from P{\"o}schel's paper
\cite{poschel:invariant-manifolds}.

Let $M$ be a uniform bound on the modulus of $F$ in a
neighbourhood of the origin.

We use multi-index notation: $l=(l_1,\,l_2)\in \mathbb{N}^2,\
|l|=l_1+l_2,\ z^l=z_1^{l_1}z_2^{l_2}$. Let $$
    F(z)=\sum_{|l|\geq 1} \overrightarrow{f_l} z^l = \Lambda z+\sum_{|l|\geq 2}\overrightarrow{f_l} z^l,
$$ where $$
    \Lambda = \begin{pmatrix} \lambda & 0 \\  0 & 1\end{pmatrix} \quad\text{and}\quad
    \overrightarrow{f_l}=(f_l^1,\,f_l^2)\in \mathbb{C}^2.
$$ The $\overrightarrow{f_l}=\overrightarrow{f_l}(\lambda)$ depend on $\lambda$; we suppress the
$\lambda$.

We shall also use the following definitions: $$
    \epsilon_n^1=\lambda^n-\lambda,\quad \epsilon_n^2=\lambda^n-1\quad\text{and}\quad \epsilon_n=\min(\epsilon_n^1,\,\epsilon_n^2)\quad\text{for}\ n\geq 2,
$$ $$
    E_n=\lambda^n Id-\Lambda=\begin{pmatrix} \epsilon_n^1 & 0 \\ 0 & \epsilon_n^2 \end{pmatrix}.
$$ For vectors in $\mathbb{R}^2$ define the usual lexicographic
ordering $$
    (z_1',\,z_2')\preceq(z_1'',\,z_2'') \text{ if and only if } z_1'\leq z_1'' \text{ and } z_2'\leq z_2''.
$$ Extend this to formal power series in $\mathbb{R}^2[[w]]$ by $$
    \sum \overrightarrow{a_i} w^i \preceq \sum \overrightarrow{b_i} w^i
    \text{ if and only if } \overrightarrow{a_i}\preceq \overrightarrow{b_i} \text{ for
    all } i.
$$ For vectors in $\mathbb{C}^2$ denote as usual
$\|z\|:=\max(|z_1|,\,|z_2|)$ and by $\{\cdot\}$ the ``norm''
$\{(z_1,\,z_2)\}:=(\|z\|,\,\|z\|)$. We view $w$ as a formal object
and extend the norms to $\mathbb{C}^2[[w]]$ by $\{\sum
\overrightarrow{a_i} w^i\}:=\sum \{\overrightarrow{a_i}\} w^i$. We
also let the usual norm on $\mathbb{C}$ act on $\mathbb{C}[[w]]$
by $|\sum \overrightarrow{a_i} w^i|:=\sum|\overrightarrow{a_i}|
w^i$. With these conventions observe that $\{A(w)+B(w)\}\preceq
\{A(w)\}+\{B(w)\}$ and $|A(w)^l|\preceq \{A(w)\}^l$ for any formal
power series $A(w),\ B(w)$ with coefficients in $\mathbb{C}^2$.
For a diagonal matrix
\[A=\begin{pmatrix} a^1 & 0 \\ 0 & a^2
\end{pmatrix}\]
define
\[\{A\}=\begin{pmatrix} \|(a^1,\,a^2)\| & 0
\\ 0 & \|(a^1,\,a^2)\|
\end{pmatrix}.\]
Observe that $\{A^{-1}\}^{-1}\{\overrightarrow{v}\}\preceq
\{A\overrightarrow{v}\}$, for any vector
$\overrightarrow{v}$.\\

We are looking for a map
\begin{equation} \label{formal psi:poeschel}
    \psi(w)=\sum_{k= 1}^\infty \overrightarrow{\psi_k} w^k=Jw+\sum_{k=2}^\infty
    \overrightarrow{\psi_k} w^k,
\end{equation}
where $$
     J=\overrightarrow{\psi_1}=(1,\,0),\quad \overrightarrow{\psi_k}=(\psi_k^1,\,\psi_k^2)
     \in \mathbb{C}^2,\quad
    k\in \mathbb{N},\quad w\in\mathbb{C},
$$ satisfying the functional equation $F\circ\psi(w)=\psi(\lambda
w)$. In power series this functional equation can be written
\begin{equation} \label{formal power ser eqn:poeschel}
    \sum_{n\geq 2}(\lambda^n Id -\Lambda) \overrightarrow{\psi_n}  w^n=\sum_{|l|\geq 2}
    \overrightarrow{f_l}\left(\sum_{k\geq1}\overrightarrow{\psi_k} w^k\right)^l.
\end{equation}
Applying  $\{\cdot\}$ to (\ref{formal power ser eqn:poeschel})
gives
\begin{align*}
    \sum_{n\geq 2} \{E_n \overrightarrow{\psi_n} \} w^n & \preceq \sum_{|l|\geq 2}
    \left\{\overrightarrow{f_l}\left(\sum_{k\geq1}\overrightarrow{\psi_k} w^k\right)^l\right\} \\
    & \preceq \sum_{|l|\geq 2} \{\overrightarrow{f_l}\} \left|\left(\sum_{k\geq 1}
    \overrightarrow{\psi_k}
    w^k\right)^l \right| \\
    & \preceq \sum_{|l|\geq 2} \{\overrightarrow{f_l}\} \left(\sum_{k\geq 1} \{\overrightarrow{\psi_k}\}
    w^k\right)^l.
\end{align*}
We conclude that
\begin{equation}  \label{formal power ser ineq:poeschel}
    \sum_{n\geq 2} \{E_n^{-1}\}^{-1} \{\overrightarrow{\psi_n}\} w^n
    \preceq \sum_{|l|\geq 2} \{\overrightarrow{f_l}\} \left(\sum_{k\geq 1} \{\overrightarrow{\psi_k}\}
    w^k\right)^l.
\end{equation}\\

Comparing powers of $w$ in (\ref{formal power ser eqn:poeschel})
it is possible to find $\overrightarrow{\psi_n}$ recursively: $$
    \overrightarrow{\psi_n}=E_n^{-1}\,P_n(\overrightarrow{\psi_1},\ldots,\,\overrightarrow{\psi_{n-1}},\,
     (\overrightarrow{f_l})_{2\leq |l|\leq n})=E_n^{-1} \,(P_n^1,\, P_n^2),\quad n\geq 2.
$$ The functions $P_n^1,\ P_n^2$ are polynomials in the
coordinates of their vector arguments. The coefficients of these
polynomials are positive. Note that the coefficients do not depend
on $\lambda$ and $\overrightarrow{f_l}$; they are the same for all
linearization problems in this context. For our purposes they are
universal polynomials.

Equating coefficients in powers of $w$ in (\ref{formal power ser
ineq:poeschel}) the same way as we did  in  (\ref{formal power ser
eqn:poeschel}) we can rewrite  (\ref{formal power ser
ineq:poeschel})  as
\begin{equation} \label{mod psi ineq:poeschel}
    \{\overrightarrow{\psi_n}\} \preceq \{E_n^{-1}\}\, P_n(\{\overrightarrow{\psi_1}\},
    \ldots,\{\overrightarrow{\psi_{n-1}}\},\, (\{\overrightarrow{f_l}\})_{2\leq |l|\leq n}),\ n\geq
    2.
\end{equation}
Define the polynomials $Q_n$ in $E_j^{-1}$ and
$\overrightarrow{f_l}$ recursively by
\begin{align*}
    Q_1  & =(1,\,0), \\
    Q_n & =Q_n(E_2^{-1}, \ldots, E_{n-1}^{-1},\,(\overrightarrow{f_l})_{2\leq |l|\leq n})\\
    & =E_n^{-1}\,P_n(Q_1, \ldots,Q_{n-1},\, (\overrightarrow{f_l})_{2\leq |l|\leq n}),\quad n\geq 2.
\end{align*}
These are the polynomials which result in ``unravelling'' the
recursive relations defining the $\overrightarrow{\psi_n}$. In
other words, they exhibit the explicit dependence of
$\overrightarrow{\psi_n}$ on the coefficients
$(\overrightarrow{f_l})_{|l|\geq2}$ of $F$ and the small divisors
$\lambda^k-\lambda,$ and $\lambda^k-1$. That is to say
$\overrightarrow{\psi_n}=Q_n$. Then from (\ref{mod psi
ineq:poeschel}) and the definition of $\{\cdot\}$ it follows that
\begin{align*}
    \|\overrightarrow{\psi_n}\| & \leq \|Q_n(\{E_2^{-1}\},\ldots,\{E_{n-1}^{-1}\},\,
    (\{\overrightarrow{f_l}\})_{2\leq |l|\leq n}) \| \\
                & \leq \|Q_n(\{E_2^{-1}\},\ldots,\{E_{n-1}^{-1}\},\, ((M,\,M))_{2\leq |l|\leq n}) \|.
\end{align*}
We thus see that in order to exhibit a uniform exponential bound
on $\overrightarrow{\psi_n}$ it is enough to exhibit such a bound
with $\overrightarrow{f_l}$ replaced by
$\overrightarrow{M}=(M,\,M)$ and the small divisors
$\epsilon_n^1,\ \epsilon_n^2$ replaced by $\epsilon_n$ in
(\ref{formal power ser ineq:poeschel}).

To do so we proceed as follows. Let
$\overrightarrow{\sigma_n}=(\sigma_n,\,\sigma_n)$ be a sequence
defined by setting $\overrightarrow{\sigma_1}=(1,\,1)$ and
equating coefficients of powers of $w$ in
\begin{equation} \label{sigma overline ineq:poeschel}
    \sum_{n\geq 2}\{(\lambda^n Id-\Lambda)^{-1}\}^{-1}\,
    \overrightarrow{\sigma_n}
    w^n=\sum_{|l|\geq 2} \overrightarrow{M} \left( \sum_{k\geq 1} \overrightarrow{\sigma_k} w^k\right)^l.
\end{equation}
In other words
\begin{align*}
    \overrightarrow{\sigma_n} & =\{E_n^{-1}\}\,P_n((\overrightarrow{\sigma_k})_{2\leq k\leq
        n-1},\,(\overrightarrow{M})_{2\leq|l|\leq n}) \\
&
=Q_n(\{E_2^{-1}\},\ldots,\{E_{n-1}^{-1}\},\,(\overrightarrow{M})_{2\leq|l|\leq
n}).
\end{align*}
 Using the fact that $P_n$ and $Q_n$ are polynomials with
positive coefficients it is clear that
$\{\overrightarrow{\psi_n}\}\preceq \overrightarrow{\sigma_n}$.
Now rewrite (\ref{sigma overline ineq:poeschel}) in terms of
$\sigma_n$:
\begin{align*}
    \sum_{n\geq 2}\epsilon_n \sigma_n w^n &=M\sum_{|l|\geq 2}  \left( \sum_{k\geq 1} \overrightarrow{\sigma_k} w^k\right)^l \\
        & = M\sum_{\nu\geq 2} \sum_{l_1+l_2=\nu}\left( \sum_{k\geq 1} \sigma_k w^k\right)^\nu \\
        & = M\sum_{\nu\geq 2} (\nu+1) \left( \sum_{k\geq 1} \sigma_k w^k\right)^\nu \\
        & = M \sum_{n\geq 2}w^n \sum_{ k_1+\cdots+k_\nu=n,\, \nu\geq 2,\ k_i\geq 1} (\nu+1)\, \sigma_{k_1}\cdots\sigma_{k_\nu}.
\end{align*}
We conclude that
\begin{equation} \label{sigma explicit recursion:poeschel}
    \sigma_n=\epsilon_n^{-1} M \sum_{\nu\geq 2,\ k_1+\cdots+k_\nu=n,\ k_i\geq 1}
    (\nu+1)\, \sigma_{k_1}\cdots\sigma_{k_\nu},\quad n\geq 2.
\end{equation}
Up until now we have been working with $\lambda=re^{i\theta_0}$
with an arbitrary fixed $r$. According to Lemma \ref{replace small
divisors lemma:poeschel} below we can replace $\lambda$ and the
corresponding small divisors $\epsilon_n$ in (\ref{sigma explicit
recursion:poeschel}) with the ones corresponding to $r=1,\
\lambda_0=e^{i\theta_0}$, by increasing $M$ by a uniform factor.
For simplicity we keep the same notation.

Following P{\"o}schel-Brjuno-Siegel (see
\cite{poschel:invariant-manifolds} page 959) we split the problem into two,
one involving no small divisors, and one involving only the small
divisors. Let $$
    \eta_1=1,\ \ \eta_n=M \sum_{\nu\geq 2,\ k_1+\cdots+k_\nu=n,\ k_i\geq 1} (\nu+1)\, \eta_{k_1}\cdots\eta_{k_\nu}
$$ and $$
    \delta_1=1,\ \ \delta_k=\frac1{\epsilon_k} \max_{\nu\geq2,\ k_1+\cdots+k_\nu=k,\ k_i\geq 1} \delta_{k_1}\cdots\delta_{k_\nu},\ k\geq 2,\ n\geq 2.
$$ It is easy to see by induction that $\sigma_n\leq\eta_n
\delta_n$. We refer to P{\"o}schel (\cite{poschel:invariant-manifolds}, pages
959-963)
for the bound $\delta_n\leq C
a^n,\ a=a(\theta_0)$, and restrict our attention to bounding
$\eta_n$. These numbers satisfy $$
    \sum_{n\geq 2}\eta_n w^n = M \sum_{\nu\geq 2} (\nu+1)\, \left( \sum_{k\geq 1} \eta_k w^k\right)^\nu.
$$ Setting $\eta=\eta(w)=\sum_{n\geq1} \eta_n w^n$ the last
equation becomes $$
    \eta-w=M\sum_{\nu\geq 2} (\nu+1)\, \eta^\nu=M\left(\frac1{(1-\eta)^2}-1-2\eta\right).
$$ By the implicit function theorem this defines $\eta=\eta(w),\
\eta(0)=0$, as an analytic function in a disk $w\in D(0,\,1/b)$.
Here $\eta=\eta(w)$ depends on $M$ only.

In particular, $\eta_n\leq C b^n$ for some $C>0$ also dependent on
$M$. Therefore $\|\psi_n\|\leq C(abM)^n$ independently of $r$.
This, together with the observation that $\psi_n$ are rational
functions in $\lambda$ and thus continuous in $r$ for an
irrational $\theta$, shows that $\psi=\psi(\cdot\,;r)$ is a
continuous family in $\mathcal{O}(D(0,\,\rho),\,\mathbb{C}^2)$,
equipped with the compact-open topology.  This $\rho=\frac1{abM}$
is the one mentioned in the statement of the theorem.

We now address the question of smoothness of this family. By the
chain rule $\frac{d}{dr}=\frac{d\lambda}{dr} \frac{d}{d\lambda}=
e^{2\pi i\theta}\frac{d}{d\lambda}$. We consider the formal
derivative of the series $\psi$, $$
    \frac{d}{d\lambda}\psi=\sum_{n\geq1}\frac{d\overrightarrow{\psi_n}}{d\lambda} w^n.
$$ To prove convergence and continuity in $r$ it is enough to
demonstrate a uniform exponential bound on the coefficients
$\frac{d}{d\lambda}\overrightarrow{\psi_n}$. To this end $$
    \frac{d}{d\lambda}\psi_n^j =\frac{d}{d\lambda}P_n^j=\frac{d}{d\lambda}Q_n^j,\ \ j=1,2.
$$ As we observed previously, $Q_n^j$ is a polynomial in
$\overrightarrow{f_l},(\epsilon_n^1)^{-1},(\epsilon_n^2)^{-1}$
with positive coefficients:
\begin{equation} \label{Qn explicit form:poeschel}
    Q_n^j=\sum A_{s,t,l,q}\, f_{l_1}\cdots f_{l_s}\ \epsilon_{q_1}^{-1}\cdots\epsilon_{q_t}^{-1}
\end{equation}
 where $(l=l_i)_{1\leq i\leq s}$ and $q=(q_j)_{1\leq j\leq t}$ are appropriate sequences of
indexes\footnote{We abuse notation here and consider $f_{l_i}$
denoting a component of $f_{l_i}$ rather than $f_{l_i}$ itself.
This is part of what we mean by ``appropriate indexes''. The other
problem that we do not address explicitly is the index sets over
which $l_i$  vary. The same disclaimer applies to the index $q_i$
of $\epsilon_{q_i}$ } and $A_{s,t,l,q}\geq 0,\ |l_i|\leq n,\
|q_i|\leq n$. We refrain from giving more details on the explicit
expansion of $Q_n$'s as this will require introducing a tree
formalism to deal with their recursive definition and is not
necessary for our purposes (\cite{chierchia-falcolini:trees}).
What is important for us is that $s,\,t\leq n^2$: none of the
monomials are of degree higher than $2n^2$. To prove this denote
by $d_n$ the maximum number of $\epsilon$'s in a monomial of
$Q_n$. Since $Q_1=(1,\,0)$ we have $d_n=0$. We proceed by
induction. Considering (\ref{formal power ser eqn:poeschel}) we
see that $d_n\leq 1+\max\,(d_{k_1}+\cdots+d_{k_\nu})$ where the
maximum extends over all $2\leq\nu\leq n$, $k_i\geq 1$,
$k_1+\cdots+k_\nu=n$. Then inductively $$
    d_n \leq 1+\max\,(k_1^2+\cdots+k_\nu^2)<1+n^2$$
and thus $d_n\leq n^2$. The last inequality follows since
$k_i\geq1$. Similarly we can prove that $t\leq n^2$.

Consider the $\frac d{d\lambda}$ derivative of (\ref{Qn explicit
form:poeschel}). The coefficients $A_{s,t,q,l}$ are independent of
$\lambda$ and the product rule gives
\begin{align*}
    \left|\frac d{d\lambda}Q_n^j\right| =  &|\sum A_{s,t,q,l}\, f_{l_1}\cdots \frac {df_{l_i}}{d\lambda}\cdots f_{l_s}\ \epsilon_{q_1}^{-1}\cdots\epsilon_{q_t}^{-1} + \\
    & + \sum A_{s,t,q,l}\, f_{l_1}\cdots  f_{l_s}\ \epsilon_{q_1}^{-1}\cdots\frac {d\epsilon_{q_k}^{-1}}{d\lambda}\cdots\epsilon_{q_t}^{-1}| \\
    \leq  & \sum A_{s,t,q,l}\, |f_{l_1}|\cdots |\frac {df_{l_i}}{d\lambda}|\cdots |f_{l_s}|\ |\epsilon_{q_1}^{-1}|\cdots|\epsilon_{q_t}^{-1}| + \\
    & + \max\,\left( \epsilon_{q_k}\frac {d\epsilon_{q_k}^{-1}}{d\lambda}\right) \sum A_{s,t,q,l}\, |f_{l_1}|\cdots |f_{l_s}|\ |\epsilon_{q_1}^{-1}\cdots\epsilon_{q_t}^{-1}|\\
    \leq & Cd_n \,|Q_n^j((M,M),\cdots, \epsilon_n^{-1},\cdots)| \\
    \leq & C d_n\,\sigma_n.
\end{align*}
The second to last inequality is easily deduced by an explicit
computation of $\epsilon_{q_k}\frac
{d\epsilon_{q_k}^{-1}}{d\lambda}$ and utilizing the Diophantine
condition on $\lambda$.

Finally, since $\sigma_n$ is bounded exponentially and uniformly
so are the coefficients of the formally derived series $\frac
d{d\lambda} \psi$. Thus convergence and continuity of
$\frac{d\psi}{dr}$ in $r$ are established. We conclude that the
family of invariant manifolds is $C^1$-smooth along the curve
$\{\arg\lambda=\theta_0,\ 1-\delta<r<1+\delta\}$.
\end{proof}

\begin{remark}
Using the same methods it is possible to show higher order
smoothness  of the family of invariant manifolds along the curve
$\{\arg\lambda=\theta_0,\, 1-\delta<r<1+\delta\}$. In fact, the
family is $C^\infty$-smooth in $r=|\lambda|$. More general
approach regions are possible too, provided a Diophantine-type
condition remains true uniformly in $r$. However, the series
$\psi$ parametrizing the invariant manifolds are definitely not
holomorphic in $\lambda$. Their coefficients are rational
functions of $\lambda$ containing all possible terms $\lambda^n-1$
in denominators and thus explode at a dense set of points on the
curve $\{|\lambda|=1\}$.
\end{remark}

\begin{lem} \label{replace small divisors lemma:poeschel}
There exist constants $c>0$ and $N\geq0$ independent of
$1-\delta<r<1+\delta$ such that for $\lambda=re^{2\pi i\theta_0}$
$$
    |\lambda^k-1|\geq\sqrt2/2|e^{2\pi i k\theta_0}-1|\geq ck^{-N}\quad\text{for all } k\geq 1.
$$
\end{lem}

\begin{proof}
It is enough to show that if $\lambda_0=e^{2\pi i\theta_0}\in S^1$
satisfies (\ref{Siegel condition}) for $c_0$ and $N_0$ then the
same estimates holds  for all $\lambda=re^{2\pi i\theta_0}$ with
$r$ close to $1$ for some $c>0$, $N\geq0$.

Let
\[I=\{k\in\mathbb{N}\,:\, |\arg(e^{2\pi i k\theta_0})|>\pi/4 \}\]

where $-\pi<\arg(e^{2\pi ik\theta_0})\leq\pi$.

For $k\in I$ we have
\[ |r^ke^{2\pi ik\theta_0}-1|>\sin(\pi/4)=\sqrt2/2.\]

For $k\in\mathbb{N}\setminus I$
\begin{align*}
    & |r^ke^{2\pi ik\theta_0}-1| \geq |e^{2\pi ik\theta_0}-1|\cos(\arg(e^{2\pi ik\theta_0}))\geq \\
    & |e^{2\pi ik\theta_0}-1|\cos(\pi/4)=\sqrt2/2|e^{2\pi ik\theta_0}-1|
    \geq \sqrt2/2ck^{-N}.
\end{align*}
\end{proof}

We now complete the proof of Theorem \ref{main theorem} The idea
of the proof is as follows: We assume that the normal eigenvalue
is not constant. Then, for certain points, $x_1$ and $x_2$, lying
in $\Sigma$ we shall create paths, $\Gamma_1$ and $\Gamma_2$. The
$\Gamma_i$ are to lie entirely in $J$, and will pass through
$x_i$. At each point on $\Gamma_i$ we will find an invariant
manifold tranverse to $J$. This family of invariant manifolds will
vary in a $C^1$ smooth fashion along $\Gamma_i$. The two families
of manifolds thus built disconnect a small neighbourhood. Using
the fact that the $\gamma$ constructed in Lemma \ref{rank 0
theorem} does not intersect the two families of manifolds we will
be able to create a recurrence in $\Omega$. This contradiction
will complete the proof.
\\

We have shown that if there is more than one limit map then
$\Sigma$ lies in $J$, and the normal eigenvalue to $J$ along
$\Sigma$ is of constant modulus $1$. We focus our attention on
smooth points of $J$ in $\Sigma$, and use the local form described
in Lemma \ref{rank 0 theorem}, where $J$ is the $x$ axis. We
assume, by way of contradiction, that the normal eigenvalue is not
constant. For ease of notation we shall call this eigenvalue
$\lambda(x,0)$ or $\lambda(x)$. We have noted before that at
smooth points of $J$, $\lambda(x)$ is holomorphic.

We separate $J$ into three sets $J_s,\, J_u,$ and $J_n$ on which
$F$ is respectively attracting, repelling or neutral in the normal
direction:
\begin{align*}
J_s & = \{(x,0)\in J\mid |\lambda(x,0)|<1\}, \\
J_u & = \{(x,0)\in J\mid |\lambda(x,0)|>1\}, \text{ and} \\
J_n & = \{(x,0)\in J\mid |\lambda(x,0)|=1\}.
\end{align*}

We now consider a generic point $(x_0,0)\in J_n$: $(x_0,0)$ is a
smooth point of $J$ where $\lambda(x_0)\neq1$ and
$\lambda'(x_0)\neq0$. By the holomorphicity of $\lambda$ and the
implicit function theorem there is a small neighborhood of
$(x_0,0)\in U\subset\mathbb{C}^2$ where
\[J_n\cap
U=\lambda^{-1}(S^1)\cap U=\{(x,0)\in J\mid|\lambda(x)|=1\}\cap U\]
is a smooth real one dimensional curve. By
\cite{khinchin:continued-fractions} inequality (\ref{Siegel
condition}) is satisfied on a full measure set $D'\subset S^1$ .
Thus $D=\lambda^{-1}(D')\subset J_n$ is dense in $J_n$, and we can
also choose $x_0\in D$, i.e. $\lambda(x_0)\in D'$, so that
(\ref{Siegel condition}) is satisfied at $(x_0,0)$ with some
$c_0>0$, $N_0\geq0$.  We consider the smooth real curve in $J$
\[ \Gamma'=\{(x,0)\in J\mid \arg\lambda(x)=\arg\lambda(x_0)\}\cap U.\]
Observe that  $\lambda'(x_0)\neq0$ implies $\Gamma'$ is
transversal to $J_n$ at $(x_0,0)$. If we choose $\Gamma$ to be a
small connected component of $\Gamma'$ containing $(x_0,0)$, then
$\Gamma$ does not intersect $J_n$, except at $(x_0,\,0)$.

By Theorem \ref{theorem:invariant manifolds} we obtain a family of
local invariant manifolds, one passing through
each point on $\Gamma$.\\

Choose an $x_0\in\Sigma$ and a neighbourhood, $U(x_0)$, of $x_0$
in $J$ satisfying the following:
\begin{enumerate}
\item $J$ is smooth in $U(x_0)$,
\item $\lambda(x)\neq 1$ in $U(x_0)$,
\item $\lambda'(x)\neq 0$ in $U(x_0)$, and
\item $J_n\cap U(x_0)$ is a smooth curve.
\end{enumerate}
This is possible since the singular points of $J$ are isolated,
and since $\lambda(x)$ is holomorphic away from singular points of
$J$ . The fact that $\lambda'(x)\neq 0$ in $U(x_0)$ means that
$\lambda$ is a local diffeomorphism in $U(x_0)$. Thus, shrinking
$U(x_0)$ if necessary, we can insure that $J_n\cap
U(x_0)=\lambda^{-1}(S^1)\cap U(x_0)$ has a single component. Let
us write $J_n\cap U(x_0)$ as $\eta:[0,\,1]\rightarrow U(x_0)$. Let
\[t_1 =\inf\{t\in [0,\,1]\mid \eta(t)\in\Sigma\},\]
and
\[t_2 =\sup\{t\in [0,\,1]\mid \eta(t)\in\Sigma\}.\]
Now pick $x_1$ and $x_2$ in $U(x_0)$ satisfying:
\begin{enumerate}
\item $x_i=\eta(s_i)$ for $i=1,\,2$ with $t_1<s_1<s_2<t_2$, and
\item $\lambda(x_i)$ satisfies the Diophantine condition (\ref{Siegel
condition}), for $i=1,\,2$.
\end{enumerate}
This is possible since Diophantine points are dense in $J_n$.

Just as we chose the path $\Gamma$ for the point $x_0$ above, we
now choose paths $\Gamma_1$ and $\Gamma_2$ for the points $x_1$
and $x_2$. We thus obtain two families of invariant manifolds, one
each along $\Gamma_1$ and $\Gamma_2$. Call these families of
manifolds $M_1$ and $M_2$ respectively. Since the manifolds vary
in a $C^1$ smooth fashion, we see that shrinking $U(x_0)$
suitably, $M_1$ and $M_2$ each disconnect $U(x_0)$. Additionally,
by further shrinking $U(x_0)$ if needed, we can assume that $M_1$
and $M_2$ are disjoint. Thus the two families and the boundary of
$U(x_0)$ bound a compact set, $V$.

We claim that the invariant curve $\gamma$ constructed in Lemma
\ref{rank 0 theorem} cannot intersect $M_1\cap U(x_0)$ or $M_2\cap
U(x_0)$. To see this we proceed as follows. For ease of notation
let $M$ represent either $M_1$ or $M_2$. Assume by way of
contradiction that $\gamma$ intersects $M\cap U(x_0)$ at the point
$p$. There are three cases:
\begin{enumerate}
\item $p$ is in a stable manifold in $M$ (a manifold
corresponding to a point in $J_s$),
\item $p$ is in a centre manifold in $M$ (a manifold
corresponding to a point in $J_n$), or
\item $p$ is in an unstable manifold in $M$
(a manifold corresponding to a point in $J_u$).
\end{enumerate}
In the first case, $F^n(p)$ converges to a single fixed point. We
have assumed, however, that $\{F^n(p)\}$ has more than one limit
point. In the second case, $p$ is a limit point of $\{F^n(p)\}$.
We have assumed, however, that $\Omega$ is nonrecurrent. In the
third case the preimages, $F^{-n}(p)$, remain in $M\cap U(x_0)$.
By shrinking $U(x_0)$ we can assume that $\gamma_0$ does not
intersect $U(x_0)$. The point $p$ however is the image $F^n(p')$
for some $p'\in\gamma_0$ and some $n>1$. Thus $p'$ must be in
$U(x_0)$, but we have assumed that it is not.

We know that $\eta([s_1,\,s_2])$ is contained in $V$, and that $V$
contains a point $q_1$ in $\Sigma$. We also know that
$\eta(t_1)\in\Sigma$ and $\eta(t_2)\in\Sigma$ are in the
complement of $V$. The invariant curve $\gamma$ comes arbitrarily
close to both of these points. Thus $\gamma$ must leave and return
to $V$ infinitely many times. Let $U'$ be a small neighbourhood of
$J_n$. Then $V\backslash \overline{U'}$ is compact, and thus
$\gamma$ has accumulation points in $V\backslash \overline{U'}$.
This is a contradiction, since all limit points of $\{F^n\}$ in
$U(x_0)$ lie on $J_n$.

We notice that if the normal eigenvalue is constant along $J$, and
satisfies the Diophantine condition (\ref{Siegel condition}), then
precisely the arguments above show that there cannot be more than
one rank $0$ limit map.

\begin{remark}
An example of an automorphism with one limit map is given in
section \ref{maps:rank 0}. At present we do not know if a there
are any maps with more than one rank $0$ limit.
\end{remark}

\subsection{One Limit Map}
We make several comments about the case where there is only one
rank $0$ limit map.

Assume there is only one rank $0$ limit map, $h$. Just as in the
case where there is more than one rank $0$ limit, $h(\Omega)$ is a
fixed point.

If $q$ is an isolated fixed point then $q$ can be attractive,
repulsive, saddle, semi repulsive, semi attractive or neutral.
However, \cite{ueda:local} Proposition 4.1 shows that the only
fixed points to which other points converge uniformly are the
attractive, semi attractive or neutral ones. We can eliminate the
possibility that the fixed point is attractive, since in this case
it would be in the interior of the Fatou component, and the Fatou
component would thus be recurrent.

If $q$ lies on a curve of fixed points we can locally write the
equation for $F$ in the form of Equation \ref{normal form 1}.
Proposition 4.1 from \cite{ueda:local} and the fact that $q$ is on
the curve of fixed points show that $q$ can only be semi
attractive or neutral. The arguments in Section \ref{structure of
sigma} show that $q$ is a neutral fixed point.

\section{The Rank $1$ Case}\label{rank 1}
We define $F$, $\Omega$ and $\Sigma$ as in Section \ref{rank 0}.
In this section we allow elements of $\Sigma$ to have rank $0$ or
$1$. We focus our attention mainly on rank $1$ elements of
$\Sigma$.

\subsection{Fibres of $h$}\label{fibers}
Let $h:=\lim_{k\rightarrow \infty} \left(F|_{\Omega}\right)^{n_k}$
be a generically rank $1$ map. We make several elementary comments
about the fibres of $h$.

\begin{definition}[$V^{h}_q$]
Let $V^{h}_{q}:=\{p\in\Omega\mid h(p)=q\}$. We may suppress the
superscript $h$ if it is clear from the context.
\end{definition}

\begin{lem} Fix $q\in h(\Omega)$. Let $V_q'$ be the pure one dimensional
irreducible components of $V_q$. Let $N_q:=\{p\in V_q'\mid
Dh(p)=0\}$. Let $V$ be an irreducible component of $V_q'$.
If $s\in N_q\cap V$ then either
\begin{enumerate}
\item $s$ is isolated in $N_q\cap V$, or
\item $V\subset N_q$.
\end{enumerate}
\end{lem}

\begin{proof}
Since $V$ is an analytic variety and $Dh$ is a holomorphic
function on $V$, the zeros of $Dh$ on $V$ are either
isolated or form an open set in $V$. In other words, $s$ is isolated
in $N_q\cap V$ or $V\subset N_q$.

\end{proof}

We also notice that at points of $p\in\Omega$ where the rank of
$Dh$ is $1$, the image of $h$ in a neighbourhood of $p$ is smooth.
This is regardless of whether $p$ lies on a fibre of $h$ which
contains points where rank $Dh$ is $0$.

Finally we note that $F(V^{h}_{q})=V^{h}_{F(q)}$, or in other
words $h\circ F = F\circ h$. This implies that $F$ is in some
sense an automorphism of $h(\Omega)$.

\subsection{Some Functional Relationships}
\begin{lem}
Let $h_1$ and $h_2$ be two rank $1$ limit maps of
$\left\{\left(F|_{\Omega}\right)^n\right\}_{n=1}^{\infty}$:
\[h_1=\lim_{k\rightarrow\infty}\left(F|_{\Omega}\right)^{m_k},\] and
\[h_2=\lim_{k\rightarrow\infty}\left(F|_{\Omega}\right)^{n_k}.\]

If $\left\{m_k-n_k\right\}$ is a finite subset of $\mathbb{Z}$,
then $h_1$ and $h_2$ have the same fibres, and in fact
$h_1=F^{l}\circ h_2$ for some $l$.
\end{lem}

\begin{proof} Assume without loss of generality that $m_k>n_k$
and pick subsequences on $\left\{m_k\right\}$ and
$\left\{n_k\right\}$ so that $m_k-n_k=l$. Then
$F^{m_k}=F^{n_k}\circ F^l$. Taking limits we see that
$h_1=F^l\circ h_2$.
\end{proof}

Given $\left\{m_k\right\}$ and $\left\{n_k\right\}$, the above
holds if we can find any subsequences of $\left\{m_{k_j}\right\}$
and $\left\{n_{k_j}\right\}$ whose differences,
$\left\{m_{k_j}-n_{k_j}\right\}$, are a finite subset of
$\mathbb{Z}$.

If we can find two subsequences of $\left\{m_k\right\}$ and
$\left\{n_k\right\}$ such that (abusing notation by not including
further subscripts)
\[m_k-n_k = l_1,\] and
\[m_j-n_j = l_2,\]
and $l_1\neq l_2$ then by the lemma we have
\[h_1=F^{l_1}\circ h_2,\] and
\[h_1=F^{l_2}\circ h_2.\]
This is turn implies that $h_2=F^{l_1-l_2}h_2$, i.e. all points in
$h_2(\Omega)$ are periodic.

\subsection{Extending and Composing Limit
Maps}\label{section:composing limit maps} In this section we
describe a formal but natural method of extending $h$ from
$\Omega$ to $h(\Omega)$ continuously along orbits.

Write $\Omega =\cup_{n=1}^{\infty}\Omega_n$ where
$\Omega_n\subset\subset \Omega_{n+1}\subset\subset\Omega$. Let
$h_1$ and $h_2$ be two limit maps of
$\left\{\left(F|_{\Omega}\right)^n\right\}_{n=1}^{\infty}$:
\[h_1=\lim_{k\rightarrow\infty}\left(F|_{\Omega}\right)^{m_k}\] and
\[h_2=\lim_{k\rightarrow\infty}\left(F|_{\Omega}\right)^{n_k}.\]
We define $h_2\circ h_1$ as follows. Fix $\left\{n_k\right\}$. For
each $k$ choose $k'=k'(k)$ so that
\begin{align*}
|F^{n_k+m_{k'}}-F^{n_k}\circ h_1|_{\Omega_k}
    &= |F^{m_{k'}}\circ
        F^{n_k}-h_1\circ F^{n_k}|_{\Omega_k}\\
    &<\frac{1}{k}.
\end{align*}
Passing to a suitable subsequence define $h_2\circ
h_1:=\lim_{k\rightarrow\infty}\left(F|_{\Omega}\right)^{n_{k_j}+m_{k'_j}}$.

\begin{lem}\label{lemma:fibres}
The map $h_2\circ h_1$ is constant on fibres of $h_1$.
\end{lem}

\begin{proof}
For ease of notation, we suppress the subscript $j$ above:
\[h_2\circ h_1=\lim_{k\rightarrow\infty}\left(F|_{\Omega}\right)^{n_{k}+m_{k'}}.\]
Let $h_1(p)=h_1(p')$. Then
\begin{align*}
|h_2 \circ h_1(p)-h_2\circ h_1(p')|
   &\leq |h_2\circ h_1(p)-F^{n_k+m_{k'}}(p)|
   +|F^{n_k+m_{k'}}(p)-F^{n_k+m_{k'}}(p')|\\
   &+|F^{n_k+m_{k'}}(p')-h_2\circ h_1(p')|.
\end{align*}

We can make the first and third summands on the right hand side of
the inequality as small as we like by choosing $k$ large, since
$\left(F|_{\Omega}\right)^{n_k+m_{k'}}\rightarrow h_2\circ h_1$ by
definition. The second summand is less than the following.
\[|h_1\circ F^{n_k}(p)-F^{n_k+m_{k'}}(p)| +
|F^{n_k+m_{k'}}(p')-h_1\circ F^{n_k}(p')| + |h_1\circ
F^{n_k}(p')-h_1\circ F^{n_k}(p)|.\] We have chosen $\Omega_k$ to
be an exhaustion of $\Omega$, so $p$ and $p'$ are in $\Omega_k$
for $k$ large. Thus

\begin{equation}\label{composition estimate}
\begin{split}
|h_1\circ F^{n_k}(p)-F^{n_k+m_{k'}}(p)|
&=|F^{n_k+m_{k'}}(p)-F^{n_k}\circ h_1(p)|\\
&\leq|F^{n_k+m_{k'}}-F^{n_k}\circ h_1|_{\Omega_k} < \frac{1}{k},
\end{split}
\end{equation}
and similarly for $p'$. We see that by choosing $k$ large we can
make the first and second summands as small as desired, and the
third summand is $0$ since $h_1(p)=h_1(p')$, and $F^{n_k}\circ
h_1=h_1\circ F^{n_k}$.
\end{proof}

\begin{lem}\label{lemma:limits}
Our definition of $h_2\circ h_1$ satisfies the following:
\[h_2\circ h_1(p)= \lim_{k\rightarrow\infty}F^{n_{k_j}}\circ h_1(p).\]
\end{lem}

\begin{proof}
We have that $h_2\circ h_1$
\[|h_2\circ h_1(p)-F^{n_{k_j}}\circ h_1(p)|\leq |h_2\circ h_1(p) -
F^{n_{k_j}+m_{k'_j}}(p)|+|F^{n_{k_j}+m_{k'_j}}(p)-F^{n_{k_j}}\circ
h_1(p)|.\] By definition of $h_2\circ h_1$ the first summand on
the right hand side can be made as small as desired by making $j$
large. Applying estimate (\ref{composition estimate}) shows that
the second summand can also be made small by choosing $j$ large.
\end{proof}

We make several notes about the above construction.
\begin{enumerate}
\item We do not know whether the maps $h_2\circ h_1$ are unique:
they might depend on the subsequence of $F^{n_k+m_{k'}}$ chosen.
\item For normality in the following setting we consider $F^n$ as
maps from $h_1(\Omega)$ to $\mathbb{C}^2$, and we also allow limit
maps to be infinite. Consider a small open set, $U$, in the
immersed variety $h_1(\Omega)$. Its preimage, $h_1^{-1}(U)=V$, is
an open set in $\Omega$. Given a subsequence $\{F^{n_k}\}$, we
pass to a convergent subsequence, giving us a limit map $h_2$,
possibly infinite. We define $h_2$ on $U$ as the restriction of
$h_2\circ h_1$ to $V\cap\Omega_n$ for some $n$. Note that this
definition is independent of the choice of $n$, by Lemma
\ref{lemma:fibres}. Using this definition and Lemma
\ref{lemma:limits} we see that $\left\{F^n\right\}$ is a normal
family on $h(\Omega)$.
\item Given any rank $1$ limit map $h$ of
$\{(F|_{\Omega})^n\}_{n=1}^{\infty}$, we can extend it to
$h(\Omega)$ by considering $h\circ h$. We can then think of the
limit map as being defined on the union of $\Omega$ and
$h(\Omega)$.
\end{enumerate}

\subsection{Limit Varieties}\label{section:limit varieties}
In this section we show that the images of two distinct rank $1$
limit maps cannot intersect except in one special case. We also
show that the image of a limit map is locally irreducible. The
images of limit maps are immersed varieties: Let $h(\Omega)$ be
the image of the limit map $h$. For every point, $p\in h(\Omega)$,
there is a neighbourhood, $U$, of $p$, and a connected component
of $U\cap h(\Omega)$ containing $p$ which is the zero set of a
holomorphic function in $U$. For ease we shall call the images of
the limit maps, limit varieties.

\begin{proof}[Proof of Theorem \ref{varieties}]
To see that two limit varieties do not intersect in $0$
dimensional sets we assume that they do. Call an intersection
point $p$. In a small neighbourhood, $W_p$, of $p$ both
$h(\Omega)$ and $g(\Omega)$ are one complex dimensional varieties,
so we can say
\[h(\Omega)\cap W_p=\{\phi^{-1}(0)\},\]
for some holomorphic function $\phi:W_p\rightarrow \mathbb{C}$.
(Here, since $g(\Omega)$ and $h(\Omega)$ are immersed varieties,
we are only examining one component of $g(\Omega)\cap W_p$ and one
component of $h(\Omega)\cap W_p$.) The restriction of $\phi$ to
$g(\Omega)\cap W_p$ is either identically zero, or has isolated
zeroes. In the former case $h(\Omega)$ and $g(\Omega)$ share an
open set. In the latter case, $\phi$ restricted to any suitably
small perturbation of $g(\Omega)$ also has zeroes. Thus if we
perturb $g(\Omega)$ into $\Omega$ we see that $h(\Omega)$ contains
points in $\Omega$. This is a contradiction.
\end{proof}

We see that by the same arguments a limit variety cannot come back
and ``hit'' itself. It is also clear that it cannot locally self
intersect: the image under $h$ of a ball in $\Omega$ must be
irreducible. Note that this does not imply that $h(\Omega)$ is a
variety.

\subsection{Fixed Points}\label{fixed points}
We consider fixed points in the limit varieties of limit maps.

\begin{lem}
Fixed points in the limit varieties of limit map cannot be
attracting or repelling in $\mathbb{C}^2$.
\end{lem}

\begin{proof}
Let $q\in h(\Omega)$ be a fixed point of $F$. For ease we let
$q=(0,0)$. $q$ cannot be an attractive fixed point; if it were it
would be interior to the Fatou component.

Assume $q$ is repelling. Let $U$ be a small neighbourhood of $q$.
Pick $p\in h^{-1}(\{q\})\backslash U$, and recall that for $j$
suitably large $F^{n_j}(p)$ is in $U$. Since under $F^{-1}$ the
origin is attractive, we see that the preimages of $F^{n_j}(p)$
remain in $U$; i.e. $p$ is not a preimage of $F^{n_j}(p)$. This
contradiction shows that $q$ is not repulsive.
\end{proof}

If $q$ lies on a curve of fixed points, then we see immediately
from Section 3.1 of \cite{ueda:local} that we can also eliminate
the possibility that $q$ is a saddle. Additionally, as mentioned
in Section \ref{structure of sigma}, if the fixed point were semi
attractive we would have normality of iterates of $F$ in a
neighbourhood of the fixed point. This is not possible since $q$
is on the boundary of $\Omega$.

\section{Polynomial Automorphisms}\label{polynomial}
We would like to know which automorphisms of $\mathbb{C}^2$ have
invariant nonrecurrent Fatou components. We show in this section
that many polynomial automorphisms will not have such Fatou
components.

The basis of our analysis is the paper by Friedland and Milnor,
\cite{friedland-milnor:automorphisms}. In this paper the authors
show that any polynomial automorphism is conjugate to
\begin{enumerate}
\item an affine map,
\item a shear, $f(z,w)=(az+p(w),\,bw+c)$ with $p$ a polynomial and $ab\neq
0$, or
\item a composition, $f_n\circ\cdots\circ f_1$, of  H{\'e}non maps,
$f_j(z,w)=(w,\, p_j(w)-a_jz)$ with $p_j$ polynomials of degree at
least two.
\end{enumerate}

A simple examination of affine maps and shears shows that their
Fatou components are either empty or all of $\mathbb{C}^2$. If the Fatou
component is all of $\mathbb{C}^2$, then it is recurrent.

H{\'e}non maps have finitely many fixed points
(\cite{friedland-milnor:automorphisms}) and thus cannot have two
rank $0$ limit maps.

By the work of Bedford and Smillie, \cite{bedford-smillie:1}, if a
H{\'e}non map is hyperbolic, then the interior of the set of
points with bounded forward orbits consists of a union of sinks.
Thus these maps cannot have invariant nonrecurrent Fatou components
with rank $1$ limit maps.

We have proved Theorem \ref{theorem:polynomial}.

\section{Examples}\label{examples}
We thank Berit Stens{\o}nes for providing invaluable help with
these examples.

We note that Weickert \cite{weickert:basins} as well as Buzzard
and Forstneric \cite{buzzard-forstneric:interpolation} have also
constructed automorphisms with prescribed jets. We not only ensure
that the automorphisms are tangent to the identity at the origin
and have the prescribed jet; we have constructed the automorphisms
carefully to ensure that they leave the $w$-axis fixed as a set.

As mentioned in the Introduction, we construct three
automorphisms: in Section \ref{maps:rank 0} we construct an
automorphism with one rank $0$ limit, in Section \ref{maps:rank 1}
we construct an automorphism with one rank $1$ limit, in Section
\ref{maps:rotation} we construct an automorphism with multiple
rank $1$ limits. Their dynamics are examined in \ref{dynamics:rank
0}, \ref{dynamics:rank 1} and \ref{dynamics:rotation}}
respectively.

\subsection{Maps}\label{maps:general} We begin with four maps:
\begin{align*}
F_1(z,w) &= (z,w+z),\\
F_2(z,w) &= (ze^{w},w),\\
F_3(z,w) &= (z,w-z),\text{ and}\\
F_4(z,w) &= (ze^{-w},w).
\end{align*}
We also introduce
\begin{align*}
b_l(z,w) &= (z,z^lw), \text{ and}\\
b_l^{-1}(z,w) &= (z,z^{-l}w),
\end{align*}
with $l\in\mathbb{Z}^{+}$. \\
Note that $b_l$ is holomorphic and one to one on
$\mathbb{C}^*\times \mathbb{C}$. It also maps $\mathbb{C}^*\times \mathbb{C}$ onto itself.\\
We let $\pi_1(z,w)=z$ and $\pi_2(z,w)=w$.\\
We see that
\[G(z,w):=F_4\circ F_3 \circ F_2 \circ F_1(z,w)=
(ze^{ze^{z+w}},w+z-ze^{z+w}).\] \\
We would like to remove all pure $z$ terms up to order $l+1$ from
$\pi_2(G(z,w))$. To do so we construct a map $F_5(z,w)=(z,w+g(z))$
with $g(z)=a_2z^2+\cdots+a_{l+1}z^{l+1}$, where we inductively
choose the $a_i$ to get rid of pure $z$ terms of degree $i$. As an
example (which we shall actually use later) we calculate this
explicitly for $l=2$.

We have that $\pi_2(G(z,w))$ is
\begin{align*}
w+z-ze^{z+w}
             &=w-z^2-\frac{z^3}{2!}-zw-z^2w-\frac{z^3w}{2!}-\frac{zw^2}{2!}
              +h_1(z,w),
\end{align*}
with
\[h_1(z,w)=-z\left(\frac{w^3}{3!}+\frac{z^3}{3!}+\frac{zw^2}{2!}+\sum_{k=4}^{\infty}\frac{(z+w)^k}{k!}\right).\]
(Note that $h_1$ does not include all of the terms which are not
pure $z$ terms. We are doing this because we will need these not
pure $z$ terms later.)\\
We have that $(\pi_1(G(z,w)))^2$ is
\begin{align*}
z^2e^{2ze^{z+w}} &= z^2 +2z^3+2z^3w+h_2(z,w),
\end{align*}
with
\[h_2(z,w)
=z^2\biggl(2z^2+2z\sum_{k=2}^{\infty}\frac{(z+w)^k}{k!}+\sum_{k=2}^{\infty}\frac{(2z)^ke^{k(z+w)}}{k!}\biggr).\]\\
Similarly $(\pi_1(G(z,w)))^3$ is
\begin{align*}
z^3e^{3ze^{z+w}} &=z^3+h_3(z,w),
\end{align*}
with
\[h_3(z,w)=z^3\biggl(\sum_{k=1}^{\infty}\frac{(3z)^ke^{k(z+w)}}{k!}\biggr).\]\\
When we add $(\pi_1(G(z,w)))^2$ to $\pi_2(G(z,w))$ we get
\[w+\frac{3z^3}{2}-zw-z^2w+\frac{3z^3w}{2}-\frac{zw^2}{2}+h_1(z,w)+h_2(z,w).\]
Now we subtract $3/2(\pi_1(G(z,w)))^3$ to $\pi_2(G(z,w))$ to get
\[w-zw-z^2w+\ord(z^4,z^3w,zw^2).\]\\

For general $l$ we see that $\pi_2(F_5\circ G(z,w))$ becomes
\[w+\ord(z^{l+2},zw).\]

In Sections \ref{maps:rank 0} through \ref{maps:rotation} we
modify the above map. We shall see in Sections \ref{dynamics:rank
0} through \ref{dynamics:rotation} that these new maps have the
desired dynamical properties.
\subsubsection{Rank $0$ Example}\label{maps:rank 0} We rewrite
$\pi_2(F_5\circ G(z,w))$ as
\[w+w\sum_{k=1}^{l+1}b_kz^k+\ord(z^{l+2},z^{l+2}w,zw^2).\]\\
We add a map $F_{5a}(z,w)=(z,we^{f(z)})$ where
$f(z)=c_1z+\cdots+c_{l+1}z^{l+1}$. We can choose the $c_k$ in such
a way as to remove all terms of the form $wz^j$ for
$j=1,\ldots,l+1$ in $F_{5a}\circ G(z,w)$.\\

We add a further map, $F_6(z,w)=(z,we^{(l+1)z})$. Let
\[\tilde{H}(z,w)=F_6\circ F_{5a}\circ F_5 \circ F_4 \circ F_3
\circ F_2 \circ F_1(z,w).\] \\
We have that
\[\tilde{H}(z,w)=\bigl(z+z^2+\ord(z^3,z^2w),\bigl(w+\ord(z^{l+2},z^{l+2}w,zw^2)\bigr)
e^{(l+1)ze^{ze^{z+w}}}\bigr).\]\\
Notice that $F_i(0,w)=(0,w)$ for $i=1,\ldots,6$ so
$\tilde{H}(0,w)=(0,w)$. Then $\tilde{H}\circ b_l(z,w):
\mathbb{C}^*\times\mathbb{C}\rightarrow
\mathbb{C}^*\times\mathbb{C}$ is well defined, one to one and
onto. Thus we also have that $b_l^{-1}\circ \tilde{H}\circ b_l:
\mathbb{C}^*\times\mathbb{C}\rightarrow\mathbb{C}^*\times\mathbb{C}$
is well defined, one to one and onto.\\
For $z\neq 0$ we have
\begin{align*}
b_l^{-1}\circ \tilde{H} \circ b_l(z,w)
 &=\bigl(ze^{ze^{z+z^lw}},\bigl(w+\ord(z^{2},z^{l+2}w,z^{l}w^2)
 \bigr)e^{(l+1)ze^{ze^{z+z^lw}}}e^{-lze^{z+z^lw}}\bigr).
\end{align*}
We notice that
$b_l^{-1}\circ\tilde{H}\circ_l(z,w)\rightarrow(0,w)$ as
$(z,w)\rightarrow (0,w)$. Thus defining
\[H(z,w)=
\begin{cases}
b_l^{-1}\circ\tilde{H}\circ{b_l}(z,w) & z\neq 0 \\
(0,w)                                 & z=0
\end{cases}
\]
we see that $H$ is an automorphism of $\mathbb{C}^2$.\\
We expand $H$:
\begin{equation}\label{h form}
H(z,w)=(z+z^2+\ord(z^3,z^{l+2}w,z^{2l+2}w^2),w+wz+\ord(z^{2},z^2w,z^lw^2)).
\end{equation}
We note for future reference some of the computations leading to
this:
\begin{align*}
ze^{ze^{z+z^lw}}&= z+z^2+\ord(z^3,z^{l+2}w,z^{2l+2}w^2).\\
ze^{z+z^lw} &=z+z^2+\ord(z^3,z^{l+1}w,z^{2l+1}w^2).
\end{align*}
Putting these together we see
\begin{align*}
e^{(l+1)ze^{ze^{z+z^lw}}}e^{-lze^{z+z^lw}}
&=1+z+\frac{3z^2}{2}+\ord(z^3,z^{l+1}w,z^{2l+1}w^2).
\end{align*}
It is now much easier to see where Equation (\ref{h form}) comes
from. As we shall see in Section \ref{maps:rank 1}, by playing a
little bit with the shears we can change the final form of $H$.

\subsubsection{Rank $1$ Example}\label{maps:rank 1} To construct
our example we use the maps $F_1$ through $F_4$. We alter $F_5$
slightly, by making it remove all pure $z^4$ terms as well as
$z^2$ and $z^3$ terms. (This was discussed in Section
\ref{maps:general}.) We leave out the map $F_{5a}$. We also leave
out $F_6$, for the moment. (We will return $F_6$ to exactly the
same place it was before; we work backwards to see where the
desired $F_6$ comes from.)
\begin{multline*}b_2^{-1}\circ F_5 \circ F_4 \circ F_3
\circ F_2 \circ F_1 \circ b_2 (z,w) =\\
\bigl(ze^{ze^{z+z^2w}},
\bigl(w-zw-z^2w+\ord(z^3,z^3w,z^3w^2)\bigr)e^{-2ze^{z+z^2w}}\bigr).
\end{multline*}
We notice the following:
\begin{align*}
(w-zw-z^2w)e^z &=w-\frac{3z^2w}{2}+\ord(z^3w).
\end{align*}
We would like to duplicate this behaviour in our maps. To do so we
add another overshear, $F_6$.
\[F_6(z,w)=(z,we^{3z}).\]
Using calculations from Section \ref{maps:rank 0} we see that
\[e^{-2ze^{z+z^2w}}e^{3ze^{ze^{z+z^2w}}}=1+z+\frac{3z^2}{2}+\ord(z^3,z^{3}w,z^{5}w^2).\]
This is not quite equal to $e^z$, but is close enough:
\begin{align*}
(w-zw-z^2w)(1&+z+\frac{3z^2}{2}+\ord(z^3,z^{3}w,z^{5}w^2))\\
 &= w-\frac{z^2w}{2}+\ord(z^3w,z^3w^2).
\end{align*}
Thus we have
\begin{align*}
b_2^{-1}\circ &F_6 \circ F_5 \circ F_4 \circ F_3 \circ F_2 \circ
F_1 \circ b_2 \,(z,w)\\
&=(ze^{ze^{z+z^2w}},w-\frac{z^2w}{2}+\ord(z^3,z^3w,z^3w^2)).
\end{align*}

\subsubsection{Rotation Examples}\label{maps:rotation}
To construct this example, we use the automorphism constructed in
Section \ref{maps:rank 1} and add a rotation. We define $H$ as
\[H(z,w)= \Theta_0 \circ b_2^{-1}\circ F_6 \circ F_5 \circ F_4 \circ F_3 \circ F_2
\circ F_1 \circ b_2\, (z,w)\]\\
where $\Theta_0(z,w)=(z,e^{i\theta_0}w)$. So we have
\[H(z,w)=(ze^{ze^{z+z^2w}},e^{i\theta_0}(w-\frac{z^2w}{2}+\ord(z^3,z^3w,z^3w^2))).\]

\subsection{Dynamics}\label{dynamics}

We note that the calculations in this section are similar to these
carried out by Weickert \cite{weickert:basins} and Hakim
\cite{hakim:transformations}.

\subsubsection{Rank $0$ Example}\label{dynamics:rank 0} We show
that $H$ has a Fatou component, $\Omega$, with the following
properties:
\begin{enumerate}
\item  $\Omega\neq\mathbb{C}^2$,
\item $\lim_{n\rightarrow \infty}H^n(p)=(0,0)$ for all
$p\in\Omega$, and
\item $(0,0)$ is in the boundary of $\Omega$.
\end{enumerate}
For ease of computation we specifically choose $l=2$. For
convenience we use the following notation:
\begin{align*}
z_n &:=\pi_1(H^n(z,w)), \text{ and}\\
w_n &:=\pi_2(H^n(z,w)).
\end{align*}
We recall the expansion of $H$:
\[(z+z^2+\ord(z^3,z^{4}w,z^{6}w^2),w+wz+\ord(z^{2},z^2w,z^2w^2)).\]
We change coordinates: $z\rightarrow \frac{-1}{z}$. In the new
coordinates $H$ becomes:
\[\left(-\frac{1}{\frac{-1}{z}+\frac{1}{z^2}+\ord\left(\frac{1}{z^3},\frac{w}{z^4},\frac{w^2}{z^6}\right)},
w-\frac{w}{z}+\ord\left(\frac{1}{z^2},\frac{w}{z^2},\frac{w^2}{z^2}\right)\right).\]
We examine $H$ on the following set:
\[U_{N,\,M}:=\left\{(z,w)\in\mathbb{C}^2\mid \re(z)>N,\,|w|<M\right\},\]
where $N$ and $M$ are large numbers. $M$ is chosen arbitrarily and
is fixed, $N$ is increased as needed in the following, though only
finitely many times. (More correctly we should write $N(M)$.) We
often write $U$ instead of $U_{N,\,M}$.

We examine the $z$ coordinate first.
\begin{align*}
z_1
&=\left(z+1+\frac{1-b}{z}+\ord\left(\frac{1}{z^2}\right)\right)
\end{align*}
with $b$ a constant, and choosing a suitably large $N$,
$\re(z)>N$. Thus by choosing $N$ suitably large we see that
\[\re(z_1)>\re(z)\]
and
\[\frac{1}{2}\leq |z_1| \leq |z|+2.\]
Assuming that $U$ is invariant we have that
\[\frac{n}{2}\leq |z_n| \leq |z|+2n.\]
Indeed, all that needs to be done to show that $U$ is invariant is
to show that $|w|$ remains less than $M$ under iteration by $H$.
We examine the $w$ coordinate:
\[w_1=w-\frac{w}{z}+\ord\left(\frac{1}{z^2},\frac{w}{z^2},\frac{w^2}{z^2}\right).\]
We rewrite this slightly
\[w_1=w\left(1-\frac{1}{z}+\ord\left(\frac{1}{z^2}\right)\right)+\ord\left(\frac{1}{z^2}\right).\]
Again by choosing $N$ large we can make
$\left(1-\frac{1}{z}+\ord\left(\frac{1}{z^2}\right)\right)$ less
than $1$, say
\[\left|\left(1-\frac{1}{z}+\ord\left(\frac{1}{z^2}\right)\right)\right|<1-\epsilon.\]
Then
\begin{align*}
|w_1|&=\left|w\left(1-\frac{1}{z}+\ord\left(\frac{1}{z^2}\right)\right)+\ord\left(\frac{1}{z^2}\right)\right|\\
     &\leq |w|(1-\epsilon)+\left|\ord\left(\frac{1}{z^2}\right)\right|.\\
\end{align*}
We note that $\left|\ord\left(\frac{1}{z^2}\right)\right|$ is
bounded on $\re(z)>N$, and decreases to $0$ as $N$ increases to
$\infty$.\\
Let
\[\alpha:=\frac{|\ord(1/z^2)|}{\epsilon}.\]
Choose $N$ so large that the following hold:
\begin{enumerate}
\item $\alpha << M$, and
\item $\alpha+\epsilon M < M$.
\end{enumerate}
Then if $|w|>\alpha$ we have that
\[|w|(1-\epsilon)+\left|\ord\left(\frac{1}{z^2}\right)\right|
\leq |w|(1-\epsilon)+|w|\epsilon \leq |w|.\]
If, on the other hand, $|w|<\alpha$ then
\begin{align*}
|w|(1-\epsilon)+\left|\ord\left(\frac{1}{z^2}\right)\right| &\leq
|w|-\epsilon|w|+\alpha\epsilon\\
&<|w|+\epsilon M \\
&<\alpha +\epsilon M\\
& <M.
\end{align*}
So we indeed have that $U$ is forward invariant under $H$.

We see that $z_n\rightarrow \infty$. Now all that remains to be
shown is that $w_n\rightarrow 0$. To do so we examine the iterates
of the $w$ coordinate slightly more carefully.
\begin{align*}
w_1&=w\left(1-\frac{1}{z}+\ord\left(\frac{1}{z^2}\right)\right)+\ord\left(\frac{1}{z^2}\right),\\
w_2&=w_1\left(1-\frac{1}{z_1}+\ord\left(\frac{1}{z_1^2}\right)\right)+\ord\left(\frac{1}{z_1^2}\right)\\
&=w\left(1-\frac{1}{z}+\ord\left(\frac{1}{z^2}\right)\right)
\left(1-\frac{1}{z_1}+\ord\left(\frac{1}{z_1^2}\right)\right)\\
&\quad+\ord\left(\frac{1}{z^2}\right)\left(1-\frac{1}{z_1}+\ord\left(\frac{1}{z_1^2}\right)\right)
+\ord\left(\frac{1}{z_1^2}\right).
\end{align*}
In general $w_n$ is
\[w_n=w\prod_{j=0}^{n-1}\left(1-\frac{1}{z_j}+\ord\left(\frac{1}{z_j^2}\right)\right)
+\sum_{j=0}^{n-1}\ord\left(\frac{1}{z_j^2}\right)\prod_{k=j+1}^{n-1}
\left(1-\frac{1}{z_k}+\ord\left(\frac{1}{z_k^2}\right)\right). \]
We examine the product
\begin{equation}\label{w product:rank 0}
\prod_{j=0}^{n-1}\left(1-\frac{1}{z_j}+\ord\left(\frac{1}{z_j^2}\right)\right).
\end{equation}
To understand convergence of this product we look at convergence
of the series
\begin{equation}\label{w sum:rank 0}
\sum_{j=0}^{n-1}\log\left(1-\frac{1}{z_j}+\ord\left(\frac{1}{z_j^2}\right)\right).
\end{equation}
But
\[\log\left(1-\frac{1}{z_j}+\ord\left(\frac{1}{z_j^2}\right)\right)
=-\frac{1}{z_j}+\ord\left(\frac{1}{z_j^2}\right).\] Choosing $N$
large we can insure that
\[\left|-\frac{1}{z_j}+\ord\left(\frac{1}{z_j^2}\right)\right|>\frac{1}{2|z_j|},\]
and thus that the sum in (\ref{w sum:rank 0}) tends to negative
infinity. Thus finally the product in (\ref{w product:rank 0})
tends to $0$.\\
We now turn our attention to the sum
\[\sum_{j=0}^{n-1}\ord\left(\frac{1}{z_j^2}\right)\prod_{k=j+1}^{n-1}
\left(1-\frac{1}{z_k}+\ord\left(\frac{1}{z_k^2}\right)\right)\]
which we rewrite as
\begin{equation}\label{alpha sum:rank 0}
\sum_{j=0}^{\infty}\alpha_j^n\ord\left(\frac{1}{z_j^2}\right)
\end{equation}
with
\[\alpha_j^n =
\begin{cases}
\prod_{k=j+1}^{n-1}\left(1-\frac{1}{z_k}+\ord\left(\frac{1}{z_k^2}\right)\right)
& j\leq n-1\\
0 & j>n-1.
\end{cases}
\]
We recall that $|z_j|\geq j/2$, so
\[\left|\ord\left(\frac{1}{z_j^2}\right)\right|\leq
\frac{C}{(j/2)^2}.\]\\
Clearly then
\[\sum_{j=0}^{\infty}\ord\left(\frac{1}{z_j^2}\right)\]
converges, and then since $|\alpha_j^n|<1$ for all $j$ and all $n$
we have that (\ref{alpha sum:rank 0}) converges as well.\\
Given $\epsilon>0$ choose $P$ so that
\[\sum_{j=P+1}^{\infty}\left|\ord\left(\frac{1}{z_j^2}\right)\right|<\frac{\epsilon}{2}.\]
Choose $n$ large enough that
\[\sum_{j=0}^{P}\left|\alpha_j^n\ord\left(\frac{1}{z_j^2}\right)\right|<\frac{\epsilon}{2}.\]
Then we have that
\[\sum_{j=0}^{\infty}\left|\alpha_j^n\ord\left(\frac{1}{z_j^2}\right)\right|<\epsilon.\]
We have proved that (\ref{alpha sum:rank 0}) is in fact $0$. This
and the fact that (\ref{w product:rank 0}) converges to $0$ show
that $w_n\rightarrow 0 $. Notice that the convergence of the
infinite sums and products above is uniform on compacts in
$U_{N,\,M}$.

We have thus succeeded in showing that $H^n(z,w)\rightarrow
(\infty,0)$, or in our original coordinates $H^n(z,w)\rightarrow
(0,0)$. This convergence is uniform on compacts, and holds at
least on the open set
\[\cup_{M>> 1}\cup_{n=0}^{\infty}H^{-n}(U_{N(M),\,M}).\]
The above set is contained in a Fatou component, $\Omega$. We
notice several things about $\Omega$.
\begin{enumerate}
\item $\Omega$ is not all of $\mathbb{C}^2$. To see this recall
that the $w$ axis is fixed by $H$. Assume $\{H^n\}_{n=1}^{\infty}$
is normal in a neighbourhood of $(0,w)$ with $w$ nonzero. In this
case it is not possible for points arbitrarily close to $(0,w)$ to
converge to $(0,0)$. But this does in fact happen, precisely to
points in $U_{N,\,M}$ which are close to $(0,w)$.
\item The full sequence $H^n$ converges uniformly on compacts in
$\Omega$ to $(0,0)$. This statement is clearly true in
$U_{N,\,M}$, and every point in $\Omega$ lands in some $U_{N,\,M}$
after sufficiently many iterations of $H$.
\end{enumerate}

\subsubsection{Rank $1$ Example}\label{dynamics:rank 1} We show
that $H$ has a Fatou component, $\Omega$, with the following
properties:
\begin{enumerate}
\item  $\Omega\neq\mathbb{C}^2$,
\item $\lim_{n\rightarrow \infty}H^n$ is a rank generic $1$ map, and
\item the $w$ axis is in the boundary of $\Omega$.
\end{enumerate}
We use the same
notation as in Section \ref{dynamics:rank 0}.\\
We recall the expansion of $H$:
\[(z+z^2+\ord(z^3,z^{4}w,z^{6}w^2),
w-\frac{z^2w}{2}+\ord(z^3,z^3w,z^3w^2)).\]\\
As before we change coordinates: $z\rightarrow\frac{-1}{z}$. In
the new coordinates $H$ becomes:
\[\left(
-\frac{1}{\frac{-1}{z}+\frac{1}{z^2}+\ord\left(\frac{1}{z^3},
\frac{w}{z^4},\frac{w^2}{z^6}\right)},w-\frac{w}{2z^2} +
\ord\left(\frac{1}{z^3},\frac{w}{z^3},\frac{w^2}{z^3}\right)
\right).\]\\
As in Section \ref{dynamics:rank 0} we examine $H$ on
\[U_{N,\,M}:=\left\{(z,w)\in\mathbb{C}^2\mid \re(z)>N,\,|w|<M\right\}.\]
Exactly the same calculations as in Section \ref{dynamics:rank 0}
reveal that $z_n\rightarrow \infty$, again contingent on $U$ being
forward invariant.\\
We examine the $w$ coordinate:
\[w_1=w-\frac{w}{2z^2}+\ord\left(\frac{1}{z^3},\frac{w}{z^3},\frac{w^2}{z^3}\right)\]
or
\[w_1=w\left(1-\frac{1}{2z^2}+\ord\left(\frac{1}{z^3}\right)\right)+\ord\left(\frac{1}{z^3}\right).\]
By repeating the arguments in Section \ref{dynamics:rank 0} (with
suitable modifications) we can see that $|w|$ remains bounded
above by $M$ under iteration by $H$. Thus $U$ is indeed forward
$H$ invariant.\\
We look more closely at the iterates of the $w$ coordinate. In
general $w_n$ is
\[w_n=w\prod_{j=0}^{n-1}\left(1-\frac{1}{2z_j^2}+\ord\left(\frac{1}{z_j^3}\right)\right)
+\sum_{j=0}^{n-1}\ord\left(\frac{1}{z_j^3}\right)\prod_{k=j+1}^{n-1}
\left(1-\frac{1}{2z_k^2}+\ord\left(\frac{1}{z_k^3}\right)\right).
\]
We examine
\begin{equation}\label{w product:rank 1}
\prod_{j=0}^{n-1}\left(1-\frac{1}{2z_j^2}+\ord\left(\frac{1}{z_j^3}\right)\right).
\end{equation}
To understand convergence of this product we examine
\[\sum_{j=0}^{n-1}\log\left(1-\frac{1}{2z_j^2}+\ord\left(\frac{1}{z_j^3}\right)\right)
=\sum_{j=0}^{n-1}\left(\frac{-1}{2z_j^2}+\ord\left(\frac{1}{z_j^3}\right)\right).\]
Choosing $N$ large insures convergence of this series. We thus see
that the product (\ref{w product:rank 1}) converges to something
finite.\\
We also see from this that the series
\begin{equation}\label{w sum:rank 1}
\sum_{j=0}^{n-1}\ord\left(\frac{1}{z_j^3}\right)\prod_{k=j+1}^{n-1}
\left(1-\frac{1}{2z_k^2}+\ord\left(\frac{1}{z_k^3}\right)\right)
\end{equation}
converges and is in fact less than
\[\sum_{j=0}^{\infty}\ord\left(\frac{1}{z_j^3}\right).\]
Again, convergence of all sums and products is uniform on compacts
in $U_{N,\,M}$.

As in Section \ref{dynamics:rank 0} we are examining a Fatou
component, $\Omega$, containing an open set
\[\cup_{M>> 1}\cup_{n=0}^{\infty}H^{-n}(U_{N(M),\,M}).\]
We have again that convergence is uniform on compacts in $\Omega$.
We must show that the limit map is in fact rank $1$. We note that
having fixed an $M$ and an $\epsilon>0$, by choosing $N$ suitably
large we can insure that the quantities (\ref{w product:rank 1})
and (\ref{w sum:rank 1}) each vary by less than $\epsilon$ in
absolute value on $U_{N(M),\,M}$. Let $h$ be the limit map
$\lim_{n\rightarrow\infty}H^n$, and let $(z,w)$ and
$(\zeta,\omega)$ be two points in $U_{N,\,M}$. For ease denote the
product
\[
\prod_{j=0}^{\infty}\left(1-\frac{1}{2z_j^2}+\ord\left(\frac{1}{z_j^3}\right)\right)
\]
corresponding to $(z,w)$ (resp. $(\zeta,\omega)$) by $P_{(z,w)}$
(resp. $P_{(\zeta,\omega)}$) and the limit
\[
\lim_{n\rightarrow\infty}\,
\sum_{j=0}^{n-1}\ord\left(\frac{1}{z_j^3}\right)\prod_{k=j+1}^{n-1}
\left(1-\frac{1}{2z_k^2}+\ord\left(\frac{1}{z_k^3}\right)\right)
\]
by $S_{(z,w)}$ (resp. $S_{(\zeta,\omega)}$). Then
\begin{align*}
|h(z,w)-h(\zeta,\omega)|&
       \geq ||wP_{(z,w)}-\omega
       P_{(\zeta,\omega)}|-|S_{(z,w)}-S_{(\zeta,\omega)}||\\
       &\geq ||(w-\omega)P_{(z,w)}|-|\omega(P_{(z,w)}-
       P_{(\zeta,\omega)})|-|S_{(z,w)}-S_{(\zeta,\omega)}||.\\
\end{align*}
Letting $\omega=0$, $z=\zeta$, and choosing $w$ very large
(perhaps enlarging $M$ and $N$) we see that the last term term is
positive. Thus $h(z,0)$ is not equal to $h(z,w)$. (Note that
$P_{(z,w)}$ is bounded away from $0$, so this is in fact
possible.) Returning to the original coordinates we see that any
point in $\Omega$ converges to a point on the $w$ axis under
iteration by $H$.

Next we show that $h(\Omega)$ is the entire $w$ axis. Let $R>0$ be
a large real number . If we restrict $|w|$ to be less than $2R$,
we can choose $z_0$ so small that (\ref{w product:rank 1}) is very
close to $1$ in modulus and (\ref{w sum:rank 1}) is very small for
all $|w|<2R$. Then for $(z_0,w)\in\Omega$ with $|w|<2R$ and $z_0$
suitably chosen, we have
\[|\pi_2(h(z_0,w))-w|<1.\]\\
Restricting our attention to
$\left\{z_0\right\}\times\left\{|w|<2R\right\}$ we think of $h$ as
a holomorphic function of one variable:
\[h:\left\{z_0\right\}\times\left\{|w|<2R\right\}\rightarrow
\left\{0\right\}\times\mathbb{C}.\]\\
Now an application of the argument principle shows that
\[\left\{0\right\}\times B(0,R)\subset
h(\left\{z_0\right\} \times \left\{|w|<2R \right\}).\]
Letting $R$
increase we see that the $w$ axis is contained in $h(\Omega)$.

Finally we show that the $w$ axis is on the boundary of the Fatou
component $\Omega$. Assume $(0,w)$ is in $\Omega$. Then we know
that high enough iterates of a small neighbourhood of $(0,w)$ are
very close to the $w$ axis:
\[H^n(B(0,\delta)\times B(w,\delta))\subset B(0,\epsilon)\times
B(\pi_2(h(0,w)),\epsilon)\]\\
for all $n$ larger than $n'$. If we change coordinates,
$z\rightarrow -1/z$, we see that for $n>n'$ the $z$ coordinate of
$H^n$ remains outside a large ball centered at the origin.\\
We now notice that the estimates we made about the growth of the
$z$ coordinate under iteration by $H$ depended on the modulus of
$z$ being large, and on the $w$ coordinate remaining bounded. In
the present setting we have met all of these requirements. Thus we
see that the real part of the $z$ coordinate grows by roughly $1$
for each iteration of $H$:
\[\re(z_1) =\re(z)+1+\ord\left(\frac{1}{z}\right).\]
We have that the $z$ coordinate remains in the complement of
$B(0,1/\epsilon)$, but we know that if this is the case then the
real part of $z$ grows by $1$ under each iteration of $H$. Thus,
after a finite number of iterates some point in the complement of
$B(0,1/\epsilon)$ moves into $B(0,1/\epsilon)$. This is a
contradiction.

\subsubsection{Rotation Examples}\label{dynamics:rotation}
The dynamics of this map is essentially the same as the dynamics
of the rank $1$ map. The difference, of course, is $\Theta_0$. We
note first that the the extra multiplicative factor
$e^{i\theta_0}$ has no effect on the estimates in Sections
\ref{dynamics:rank 0} and \ref{dynamics:rank 1}. The rotation,
however, does make the full sequence, $\{H^n\}$, not convergent.

If $\theta_0$ is a rational multiple of $2\pi$ we obtain finitely
many limit maps of $H^n$. Each such map has as its image the $w$
axis. $H$ acts as periodic rotation on the $w$ axis, and the limit
maps differ from each other by composition with $H^j,$ for some
$j$.

If $\theta_0$ is an irrational multiple of $2\pi$ we obtain
infinitely many limit maps of $H^n$. Each such map has as its
image the $w$ axis. $H$ acts as an irrational rotation on the $w$
axis, and the limit maps differ from each other by composition
with $H^j,$ for some $j$.

\nocite{poschel:invariant-manifolds}
\nocite{bedford-smillie:external}
\nocite{chierchia-falcolini:trees}
\bibliographystyle{amsalpha}
\bibliography{bibliography}

\providecommand{\bysame}{\leavevmode\hbox to3em{\hrulefill}\thinspace}
\providecommand{\MR}{\relax\ifhmode\unskip\space\fi MR }
\providecommand{\MRhref}[2]{%
  \href{http://www.ams.org/mathscinet-getitem?mr=#1}{#2}
}
\providecommand{\href}[2]{#2}
\begin{thebibliography}{BS91b}

\bibitem[BF00]{buzzard-forstneric:interpolation}
Gregery~T. Buzzard and Franc Forstneric, \emph{{An interpolation theorem for
  holomorphic automorphisms of $\mathbb{C}^n$}}, The Journal of Geometric
  Analysis \textbf{10} (2000), 101--108.

\bibitem[BS91a]{bedford-smillie:1}
Eric Bedford and John Smillie, \emph{{Polynomial diffeomorphisms of
  $\mathbb{C}^2$: currents, equilibrium measure and hyperbolicity}},
  Inventiones Mathematicae \textbf{103} (1991), 69--99.

\bibitem[BS91b]{bedford-smillie:2}
\bysame, \emph{{Polynomial diffeomorphisms of $\mathbb{C}^2$. II: stable
  manifolds and recurrence}}, Journal of the American Mathematical Society
  \textbf{4} (1991), 657--679.

\bibitem[BS99]{bedford-smillie:external}
\bysame, \emph{{External rays in the dynamics of polynomial automorphisms of
  $\mathbb{C}^2$}}, {Complex geometric analysis in Pohang (1997)}
  \textbf{Contemporary Mathematics} (1999), no.~222, 41--79.

\bibitem[CF94]{chierchia-falcolini:trees}
L.~Chierchia and C.~Falcolini, \emph{{A Direct Proof of a Theorem by Kolmogorov
  in Hamiltonian Systems}}, Annali Della Scuola Normale Superiore di Pisa
  \textbf{21} (1994), 541--593.

\bibitem[CG93]{carleson-gamelin:dynamics}
Lennart Carleson and Theodore Gamelin, \emph{Complex dynamics},
  Springer-Verlag, New York, 1993.

\bibitem[FM86]{friedland-milnor:automorphisms}
S.~Friedland and J.~Milnor, \emph{{Dynamical properties of plane polynomial
  automorphisms}}, Ergodic Theory and Dynamical Systems \textbf{9} (1986),
  67--99.

\bibitem[FS95]{fornaess-sibony:fatou2}
John~Erik Forn{\ae}ss and Nessin Sibony, \emph{{Classification of recurrent
  domains for some holomorphic maps}}, Mathematische Annalen \textbf{301}
  (1995), 813--820.

\bibitem[FS98]{fornaess-sibony:fatou}
\bysame, \emph{{Fatou and Julia sets for entire mappings in $\mathbb{C}^k$}},
  Mathematische Annalen \textbf{311} (1998), 27--40.

\bibitem[Hak98]{hakim:transformations}
Monique Hakim, \emph{{Analytic transformations of $(\mathbb{C}^p,0)$ tangent to
  the identity}}, Duke Mathematical Journal \textbf{92} (1998), 403--428.

\bibitem[Khi64]{khinchin:continued-fractions}
Alexander Khinchin, \emph{{Continued Fractions}}, Phoenix Books, Chicago, 1964.

\bibitem[Nis83]{nishimura:automorphisms}
Yasuichiro Nishimura, \emph{{Automorphismes analytiques admettant des
  sous-vari{\'{e}}t{\'{e}}s de points fix{\'{e}}s attractive dans la direction
  transversale}}, Journal of Mathematics of Kyoto University \textbf{23}
  (1983), 289--299.

\bibitem[P{\"o}s86]{poschel:invariant-manifolds}
J{\"u}rgen P{\"o}schel, \emph{{On invariant manifolds of complex analytic
  mappings near fixed points}}, Expositiones Mathematicae \textbf{4} (1986),
  97--109.

\bibitem[Ued86]{ueda:local}
Tetsuo Ueda, \emph{{Local structure of analytic transformations of two complex
  variables}}, Journal of Mathematics of Kyoto University \textbf{26} (1986),
  233--261.

\bibitem[Ued94]{ueda:fatou}
\bysame, \emph{{Fatou sets in complex dynamics on projective space}}, Journal
  of the Mathematical Society of Japan \textbf{46} (1994), 545--555.

\bibitem[Wei98]{weickert:basins}
Brendan Weickert, \emph{{Attracting basins for automorphisms of
  $\mathbb{C}^2$}}, Iventiones Mathematicae \textbf{132} (1998), 581--605.

\end{thebibliography}

\end{document}